\pdfoutput=1
\documentclass[11pt,leqno]{article}

\usepackage{amssymb}
\usepackage{amsmath}
\usepackage{sectsty}
\usepackage{amsthm}
\usepackage{graphicx}
\usepackage{caption}
\usepackage{float}

\usepackage{hyperref}
\allsectionsfont{\sf}

\newif\ifTrans
\Transtrue

\usepackage{perpage} 
\MakePerPage{footnote} 

\title{\sf 
\ifTrans{Memoir on the vibratory movement\\ of an elliptical membrane\footnote{Translation of \'{E}mile Mathieu's ``M\'{e}moire sur le mouvement vibratoire d'une membrane de forme elliptique'' (\url{http://sites.mathdoc.fr/JMPA/PDF/JMPA_1868_2_13_A8_0.pdf}) by Robert H. C. Moir. Edited by Robert M. Corless.}}
\else{M\'{e}moire sur le mouvement vibratoire\\ d'une membrane de forme elliptique}
\fi}
\author{
\sf 
\ifTrans{By Mr.\ }\else{Par M.\ }\fi\'{E}mile Mathieu
}	
\date{
\sf\footnotesize 
\ifTrans{This Memoir was exhibited in January 1868 in a course at the Sorbonne}
\else{Ce M\'{e}moire a \'{e}t\'{e} expos\'{e} au mois de janvier 1868 dans un cours \`{a} la Sorbonne}
\fi
}

\begin{document}

\renewcommand*{\thefootnote}{\fnsymbol{footnote}}

\maketitle

\ifTrans
Imagine a membrane stretched equally in all directions, and
whose outline, invariably fixed, is an ellipse. 
Our goal, 
in this Memoir,   
is to determine by analysis all the circumstances
of its vibratory movement; we calculate the shape and the position
of nodal lines and the corresponding sound. But these movements are
subject to certain general laws which may be defined without
the help of analysis.

When the elliptical membrane is vibrated, there are
two systems of nodal lines which are, one of 
ellipses, the
other of 
hyperbolas, and all these curves of second order have the
same foci as the ellipse of the contour.

All these vibratory movements can be divided in two
kinds. In one of these kinds, the major axis remains fixed and forms a
nodal line, and if we consider two symmetrical points with respect to
the major axis, their movements are equal and in opposite directions. In
the other kind, on the contrary, the ends of the major axis located between
the foci and the vertices form bellies [\emph{ventres}] of vibration, while
the part located between the two foci offers a minimum of vibration,
\fi
\ifTrans
so that if we take a point $M$ on the 
line segment 
of the foci, and a
very close point on a perpendicular at $M$, the amplitude of the 
vibration is less for the first than for the second point; if we
consider any two points of the membrane, symmetrical 
relative to the major axis, their movements are equal and the same
direction.

Let us define a hyperbolic line as the two branches of a hyperbola
terminated at the major axis which have the same asymptote, so
that a hyperbola counts as two hyperbolic lines;
but if one of the axes of the membrane is stationary, it will be counted
for a single hyperbolic nodal line. So the movements of the
two kinds can be grouped two by two in a very remarkable way. 
Indeed, to a number $a$ of elliptical nodal lines and to
a number $b$ of hyperbolic nodal lines corresponds a vibrational
state of each kind. Now, although these vibrational states differ
both by the two systems of nodal lines and by the resulting sound,
they nevertheless merge into the circular membrane to
give, as lines of nodes, $a$ concentric circles and $b$ diameters, 
which divide them into equal parts. We understand from this
that if the eccentricity is very small, the sounds of these two 
vibratory states will differ very little.

We must put aside the case where there are no hyperbolic
nodal lines; because the movement cannot be of the second
kind, and there is only one vibrational state that produces $a$ nodal ellipses.

The vibratory movement of a membrane enclosed between two
confocal ellipses, all of whose points are perfectly fixed, is
also subject to very simple laws.

The nodal lines of this membrane are still ellipses and
portions of hyperbola branches which have the same foci as the
two ellipses of the contours. And there are still two kinds of
vibratory movements: in one, the portions of the major axis enclosed
between the two contours are nodes; in the other, bellies of
vibration. But when we study the vibrational states of the two kinds
which give for nodes $a$ ellipses and $b$ hyperbolic lines, we
find, if the number $b$ is large enough and the eccentricity is not very
large, that the sound is almost the same, as well as the arrangement
\fi
\ifTrans
of the nodal ellipses. Now, the two sounds differing excessively little, we
know 
that in experience we will produce the two vibrarational
states together, and, in the resulting movement, the arrangement of the $b$
hyperbolic nodal lines can vary in infinite ways.

Mr. Bourget gave the theory of the circular membrane (Annals of the
\'{e}cole Normale, t. III) and did the experiments necessary to verify 
it; he found sounds a little higher than indicated by the calculation.

\begin{center}
-----
\end{center}

{\bf 1}. Let us consider a flat, homogeneous membrane, equally stretched
in all directions, and whose contour is invariably fixed. Let us trace
in the plane of this membrane two axes of arbitrary rectangular 
coordinates, $Ox$ and $Oy$, and let us run a $z$-axis perpendicular
to this plane. If we communicate a vibratory movement to this
membrane, a point on its surface whose coordinates are $x$, $y$ and
$z = 0$ will experience a normal displacement $w$ governed by the equation
\fi
\begin{equation}
\tag{$a$}\label{eq:1a}
\tfrac{d^2w}{dt^2} = m^2\left(\tfrac{d^2w}{dx^2} + \tfrac{d^2w}{dy^2}\right),
\end{equation}
\ifTrans
where $m^2$ denotes the ratio of the tension to the density of the membrane.[\footnote{Theory of Elasticity of Mr. Lam\'{e}, IXth Lesson.}]
And we have to integrate this equation, by imposing the condition that $w$ is
zero on the contour.

We must assume in this Memoir that this outline is an
ellipse. But we will first take it to be circular and present
very succinctly the solution of this particular case, which will 
sometimes be useful to us as a means of comparison.

\begin{center}
\emph{Circular membrane.}
\end{center}

{\bf 2}. Place the origin of the coordinates in the center of the circle and pass
from the rectilinear coordinates $x$ and $y$ to polar coordinates $r$ and $\alpha$ by
\fi
\ifTrans
the formulas
\fi
$$
x = r\cos\alpha,\quad y = r\sin\alpha,
$$
\ifTrans
arbitrarily taking the direction of the polar axis.

Equation (a) becomes
\fi
$$
\tfrac{d^2w}{dt^2} = m^2\left(\tfrac{d^2w}{dr^2} + \tfrac{1}{r}\tfrac{dw}{dr} + \tfrac{1}{r^2}\tfrac{d^2w}{d\alpha^2}\right),
$$
\ifTrans
and if we set $w = u\sin 2\lambda mt$, we have
\fi
$$
\tfrac{d^2u}{dr^2} + \tfrac{1}{r}\tfrac{du}{dr} + \tfrac{1}{r^2}\tfrac{d^2u}{d\alpha^2} = -4\lambda^2u.
$$
\ifTrans
Let $u = PQ$, then by designating by $P$ a function of $\alpha$ and by $Q$ a
function of $r$ we will have an equation which can be written
\fi
$$
\tfrac{r^2}{Q}\tfrac{d^2Q}{dr^2} + \tfrac{r}{Q}\tfrac{dQ}{dr} + 4\lambda^2r^2 = -\tfrac{1}{P}\tfrac{d^2P}{d\alpha^2};
$$
\ifTrans
as the first member depends only on $r$ and the second only on $\alpha$,
they are equal to the same constant $n^2$, and we have
\fi
\begin{equation}
\tfrac{d^2P}{d\alpha^2} + n^2P = 0, \tag{1}\label{eq:11}
\end{equation}
\begin{equation}
r^2\tfrac{d^2Q}{dr^2} + r\tfrac{dQ}{dr} - \left(n^2 - 4\lambda^2r^2\right)Q = 0. \tag{2}\label{eq:12}
\end{equation}

\ifTrans
We have taken for the constant a positive quantity, in order
to obtain for $P$ the periodic function
\fi
$$
P = A\cos n\alpha + B\sin n\alpha,
$$
\ifTrans
and, so that $P$ does not change when we replace $\alpha$ by $\alpha + 2\pi$,
$n$ must be an integer.

If we integrate equation \eqref{eq:12} by series, we obtain the two particular
\fi
\ifTrans
solutions:
\fi
\begin{equation}
\tag{3}\label{eq:13}
\left\{
\begin{array}{rcl}
         Q &=& Cr^{n}\left[1 - \tfrac{(\lambda r)^2}{1\,(n+1)}
                             + \tfrac{(\lambda r)^4}{1\cdot 2\,(n+1)(n+2)}\right.\\
             &&\qquad\qquad \left.-\tfrac{(\lambda r)^6}{1\cdot 2\cdot 3\,(n+1)(n+2)(n+3)} + \ldots\right],
\end{array}
\right.
\end{equation}
\begin{equation}
\tag{4}\label{eq:14}
\left\{
\begin{array}{rcl}
         Q' &=& C'r^{-n}\left[1 + \tfrac{(\lambda r)^2}{1\,(n-1)}
                             + \tfrac{(\lambda r)^4}{1\cdot 2\,(n-1)(n-2)}\right.\\
             &&\qquad\quad\ \ \left.+\tfrac{(\lambda r)^6}{1\cdot 2\cdot 3\,(n-1)(n-2)(n-3)} + \ldots\right],
\end{array}
\right. 
\end{equation}
\ifTrans
the second of which is deduced from the first by the change of $n$
into $-n$. If we sum them up, we get the general solution; But
as obviously the vibratory movement must remain finite in the center
of the circle, and that $Q'$ becomes infinite for $r = 0$, we must confine ourselves to take for $Q$ the first particular solution, which we will carry into
\fi
\begin{equation}
\tag{5}\label{eq:15}
u = PQ,\quad w = u\sin2\lambda mt.
\end{equation}

\ifTrans
Finally, for $w$ to be a possible solution, $Q$ must be zero
along the contour circle $r = h$, and $\lambda$ is determined by the equation
\fi
$$
1 - \tfrac{(\lambda h)^2}{1\,(n+1)}
  + \tfrac{(\lambda h)^4}{1\cdot 2\,(n+1)(n+2)}
  - \cdots = 0.
$$

\ifTrans
Let us put down this equation
\fi
\begin{equation}
\tag{6}\label{eq:16}
1 - \tfrac{\tau^2}{1\,(n+1)}
  + \tfrac{\tau^4}{1\cdot 2\,(n+1)(n+2)}
  -\tfrac{\tau^6}{1\cdot 2\cdot 3\,(n+1)(n+2)(n+3)}
  + \cdots = 0;
\end{equation}
\ifTrans
it has an infinity of roots $\tau_1$, $\tau_2$, $\tau_3$,\ldots, which we will suppose 
to be arranged in increasing order of magnitude, and $\lambda$ can take any 
of the values
\fi
$$
\lambda_1=\tfrac{\tau_1}{h},\quad \lambda_2=\tfrac{\tau_2}{h},\quad\lambda_3=\tfrac{\tau_3}{h}, \ldots.
$$

\ifTrans
Thus formula \eqref{eq:15} represents an infinity of vibrational movements, 
which depend on $n$ and $\lambda$; $n$ is susceptible to all
integer values, and to each value of $n$ corresponds an infinity
of values of $\lambda$.
\fi

\ifTrans
{\bf 3}. Consider one of these vibrational states and see what the 
nodal lines are. The vibrational movement satisfies the equation
\fi
\begin{equation}
\tag{$b$}\label{eq:1b}
w = (A\cos n\alpha + B\sin n\alpha)\,Q\sin 2\lambda mt;
\end{equation}
\ifTrans
$n$ and $\lambda$ are known, and the pitch, or the number of vibrations
that occur per unit time, is $N = \tfrac{\lambda m}{\pi}$. To get the
rows of nodes, we will let $w = 0$, which we can satisfy, whatever
let t, by letting
\fi
\begin{equation}
\tag{7}\label{eq:17}
\qquad A\cos n\alpha + B\sin n\alpha = 0,
\end{equation}
\ifTrans
or by letting
\fi
\begin{equation}
\tag{8}\label{eq:18}
Q = 0.
\end{equation}
\ifTrans
From equation \eqref{eq:17} we draw $\tan n\alpha = -\tfrac{A}{B}$: therefore if we designate by $n\alpha$,
the smallest of the arcs whose tangent is $-\tfrac{A}{B}$ and by $k$
any integer, $w$ is zero for
\fi
$$
\alpha = \alpha_1 + \tfrac{k\pi}{n},
$$
\ifTrans
and therefore we have for nodal lines $n$ diameters which divide the
circumference of the circle into equal parts.

Moving on to equation \eqref{eq:18}, we first notice that $Q$ is zero at the
center of the membrane, unless $n$ is zero because of the factor $r''$,
and then it is zero for different values of $r$ which are
\fi
$$
r = \tfrac{\tau_1}{\lambda},\ \tfrac{\tau_2}{\lambda},\ \tfrac{\tau_3}{\lambda},\ldots;
$$
\ifTrans
these are the rays of the nodal circles which have the same center as the membrane.

The number of these values of $r$ for which $Q$ is zero is infinite;
but we must reject all those that are larger than the radius
of the membrane. When in formula \eqref{eq:1b} we gave ourselves the value
\fi
\ifTrans
of $n$, $\lambda$ is susceptible of an infinity of values $\tfrac{\tau_1}{h}, \tfrac{\tau_2}{h},\ldots$; 
suppose that the one we adopted is the $s^{th}$,
\fi
$$
\lambda_s = \tfrac{\tau_s}{h};
$$
\ifTrans
then the nodal circles with the number of $s - 1$ will have as radii
\fi
$$
\tfrac{\tau_1}{\lambda_s},\quad \tfrac{\tau_2}{\lambda_s},\ldots,\quad\tfrac{\tau_{s-1}}{\lambda_s}.
$$

\ifTrans
We see from the above that in examining these vibrational movements
there are no other difficulties in calculation than finding the
roots of equation \eqref{eq:1b}, in which the whole number $n$ varies.
Mr. Bourget gave, in his Memoir, a method for easily calculating
the roots of this equation, and he gave the numerical values
of these first roots for $n = 0, 1, 2, ..., 7$.

The vibrational movements represented by formula \eqref{eq:1b} are called
\emph{simple} movements, and these are the ones that are observed by
experience. Finally any vibratory movement that one can imagine is the
superposition of a finite or infinite number of these simple movements.

\begin{center}
\emph{Passage from rectilinear coordinates to coordinates of the ellipse.}
\end{center}

{\bf 4}. Denote by $A$ the semi-major axis of the elliptical membrane,
and by $c$ the half-distance between 
the foci; take as axes of $x$ and $y$
the axes of symmetry of the ellipse; then adopt a second system of
coordinates determined by the ellipses and hyperbolas which have the
same foci as the contour of the membrane.

Any one of these ellipses is given by the equation
\fi
\begin{equation}
\tag{1}\label{eq:21}
\tfrac{x^2}{\rho^2}+\tfrac{y^2}{\rho^2-c^2} = 1,
\end{equation}
\ifTrans
in which $\rho$ is $> c$, and if we let
\fi
$$
\rho=c\,\tfrac{e^{\beta}\,+\,e^{-\beta}}{2},\quad \rho'=\sqrt{\rho^2-c^2}=c\,\tfrac{e^{\beta}\,-\,e^{-\beta}}{2},
$$
\ifTrans
$\rho$ and $\rho'$ are the semi-major axis and the semi-minor axis of this ellipse, and $\beta$
is what Mr. Lam\'{e} calls the \emph{thermometric parameter} (\emph{On inverse functions of the transcendent}, I$^{\rm{st}}$ Lesson).

Any of the confocal hyperbolas has the equation
\fi
\begin{equation}
\tag{2}\label{eq:22}
\tfrac{x^2}{\nu^2}-\tfrac{y^2}{c^2-\nu^2} = 1,
\end{equation}
\ifTrans
where $\nu$ is $< c$, and if we let
\fi
$$
\nu = c\cos\alpha,\quad\nu'=\sqrt{c^2-\nu^2}=c\sin\alpha:
$$
\ifTrans
$\nu$ and $\nu'$ are the half-axes of this hyperbola, and $\alpha$ its thermometric 
parameter.

We pass from the coordinates $x$ and $y$ to the coordinates $\nu$ and $\rho$ or 
$\alpha$ and $\beta$ by means of the formulas
\fi
\begin{equation}
\tag{3}\label{eq:23}
\left\{\begin{array}{p{0.1em}p{0.05em}p{0.8em}p{0.5em}p{0.2em}}
               $x$&$=$&$\,\tfrac{\rho\nu}{c}$&$=$&$c\,\tfrac{e^{\beta}\,+\,e^{-\beta}}{2}\cos\alpha,$\\
               $y$&$=$&$\tfrac{\rho'\nu'}{c}$&$=$&$c\,\tfrac{e^{\beta}\,-\,e^{-\beta}}{2}\sin\alpha,$
 \end{array}\right.
\end{equation}
\ifTrans
that we deduce from equations (1) and (2). If we wanted to have formulas that could apply immediately to the circle, we would adopt
\fi
\begin{equation}
\tag{4}\label{eq:24}
x = \rho\cos\alpha,\quad y = \rho'\sin\alpha.
\end{equation}

\ifTrans
Let $M$ be a point which comes from the intersection of the ellipse $\beta = \beta_1$, and
the hyperbola $\nu = \nu_1$. Extend the ordinate of point $M$ until it 
meets in $N$ with the circle described on the major axis. We see from 
equations (4) that the angle $\alpha$ is equal to the angle made by the ray led 
from center to point $N$ with the $x$-axis, and as this angle has for
its cosine $\tfrac{\nu}{c}$ it is also the one made with the $x$-axis by the asymptote to the
hyperbola branch which contains the point $M$, and the ray led from the
center to the point $N$ is this asymptote.

It also follows from formulas \eqref{eq:24}, that we will obtain all the points
of the plane by supposing $\rho$ and $\rho'$ positive, and causing $\rho'$ to vary from $0$ to $\infty$, and $\alpha$ from $0$ to $2\pi$.
\fi

\ifTrans
When the coordinates are thus varied, the equation $\beta = \rm{const}$.
represents an entire ellipse, but $\alpha = \rm{const}$. does not represent any more than
one of the four branches of the hyperbola ended at the transverse axis, and the whole hyperbola is given by the four equations
\fi
$$
\alpha = \alpha_1,\quad \alpha = \pi - \alpha_1,\quad \alpha = \pi + \alpha_1,\quad  \alpha = 2\pi - \alpha_1,
$$
\ifTrans
which are those of the four branches. We assume in what follows that $\beta$ is positive; however not only is this assumption not essential, but we will have occasion to recognize in the sequel that it can be useful to vary the sign of this coordinate.

It is also good to consider the limit positions of these ellipses and these hyperbolas; for $\beta = 0$, the ellipse is reduced to the line segment 
which joins the foci $F$ and $F'$; the equation $\alpha = 0$ represents the line $Fx$ bounded at $F$ and unbounded 
in the direction of positive $x$, $\alpha = \pi$ represents the line $F'x'$ unbounded
in the direction of negative $x$; finally $\alpha = \tfrac{\pi}{2}$ determines the entire positive $y$ axis, and $\alpha = \tfrac{3\pi}{2}$ the negative part of the $y$ axis.

{\bf 5}. Let us take the equation
\fi
\begin{equation}
\tag{5}\label{eq:25}
m^2\left(\tfrac{d^2w}{dx^2}+\tfrac{d^2w}{dy^2}\right) = \tfrac{d^2w}{dt^2},
\end{equation}
\ifTrans
which by the substitution of
\fi
$$
w = u\sin 2\lambda mt
$$
\ifTrans
becomes
\fi
\begin{equation}
\tag{6}\label{eq:26}
\tfrac{d^2u}{dx^2}+\tfrac{d^2u}{dy^2} = -4\lambda^2u,
\end{equation}
\ifTrans
and substitute for $x$ and $y$ the coordinates $\alpha$ and $\beta$.

To simplify, let
\fi
$$
E(\beta) = \tfrac{e^{\beta}\,+\,e^{-\beta}}{2},\quad \mathcal{E}(\beta) = \tfrac{e^{\beta}\,-\,e^{-\beta}}{2}
$$
\ifTrans
and
\fi
$$
H^2 = E^2(\beta)\sin^2\alpha + \mathcal{E}^2(\beta)\cos^2\alpha = E^2(\beta)-\cos^2\alpha.
$$
\ifTrans
we have (II$^{\rm{nd}}$ Lesson of curvilinear coordinates of Mr. Lam\'{e})
\fi
$$
\tfrac{d^2u}{dx^2}+\tfrac{d^2u}{dy^2}=
        \left[\big(\tfrac{d\alpha}{dx}\big)^2+
              \big(\tfrac{d\alpha}{dy}\big)^2\right]
        \left(\tfrac{d^2u}{d\alpha^2}+\tfrac{d^u}{d\beta^2}\right)
$$
\ifTrans
and by differentiating the equations \eqref{eq:23}, we find
\fi
$$
\tfrac{d\beta}{dx}=-\tfrac{d\alpha}{dy}=\tfrac{\mathcal{E}(\beta)\cos\alpha}{cH^2},
\quad
\tfrac{d\alpha}{dx}=\tfrac{d\beta}{dy}=\tfrac{-E(\beta)\sin\alpha}{cH^2},
$$
$$
\big(\tfrac{d\alpha}{dx}\big)^2+ \big(\tfrac{d\alpha}{dy}\big)^2 =
\tfrac{1}{c^2H^2},
$$
\ifTrans
and we have, instead of equation \eqref{eq:26},
\fi
$$
\tfrac{1}{c^2H^2}\left(\tfrac{d^2u}{d\alpha^2}+\tfrac{d^2u}{d\beta^2}\right)=
-4\lambda^2u.
$$

\ifTrans
Let
\fi
$$u = PQ,$$
\ifTrans
and regard 
$P$ as a function of $\alpha$ only 
and $Q$ as a function
of $\beta$ only, 
and we will have, instead of the previous equation,
\fi
$$
\tfrac{d^2P}{d\alpha^2}Q+P\tfrac{d^2Q}{d\beta^2} =
-4\lambda^2c^2\left[E^2(\beta)-\cos^2\alpha\right]
$$
\ifTrans
or
\fi
$$
-\tfrac{1}{P}\tfrac{d^2P}{d\alpha^2}+4\lambda^2c^2\cos^2\alpha =
+\tfrac{1}{Q}\tfrac{d^2Q}{d\beta^2}+4\lambda^2c^2E^2(\beta).
$$

\ifTrans
Since the first member can only contain $\alpha$, and the second only $\beta$,
they are equal to the same constant $N$; so that we have, instead of a
partial difference equation, two second order differential equations
\fi
$$
\tfrac{d^2P}{d\alpha^2}+\left(N-4\lambda^2c^2\cos^2\alpha\right)P = 0,
$$
$$
\tfrac{d^2Q}{d\beta^2}-\left[N-4\lambda^2c^2E^2(\beta)\right]Q = 0,
$$

\ifTrans
The first of these equations is suitable for the circular membrane, 
if we make $c=0$, and we know that we must then take for the constant $N$
the square of an integer; it does not follow that
\fi
\ifTrans
the same thing takes place here; because we do not see that 
the equation
does not
depend on $\lambda c$; but we are sure that if the constant depends on
this quantity, it reduces at least to the square of an integer for $c = 0$.

Suppose that we know one of the values of $N$, and that we have found
values of $P$ and $Q$, which satisfy these two equations; then the
formula
\fi
$$
w= PQ\sin 2\lambda mt
$$
\ifTrans
represents a possible vibratory movement of the membrane, if we
determine $\lambda$ by the condition that $Q$ is zero for the value of
$\beta$ relative to the contour.
\fi

\ifTrans
\begin{center}
\emph{On the determination of the constant $N$.}
\end{center}

{\bf 6}. The first question we must ask ourselves 
is therefore to determine
the constant $N$. Now, for the expression of $w$ to be admitted, it is
necessary that when we change $\alpha$ to $\alpha+2\pi$, $w$ remains the
same, since $w$ will continue to give the displacement of the same point
of the membrane; thus $P$ is necessarily a periodic function whose period
is $2\pi$, and this condition must determine the constant $N$.

Let us recall the results obtained by Sturm on second order linear differential equations. Any such equation can be put in the form
\fi
\begin{equation}
\tag{1}\label{eq:31}
\tfrac{d\left(L\tfrac{dy}{dx}\right)}{dx}+Gy = 0,
\end{equation}
\ifTrans
$L$ and $G$ being two functions of $x$. Let us imagine that $G$ also
contains a parameter $h$, and give it an increment $\delta h$; then $G$
will take the value $G+\delta G$, and $y$ will change to the function $y_1$,
which satisfies the equation
\fi
\begin{equation}
\tag{2}\label{eq:32}
\tfrac{d\left(L\tfrac{dy_1}{dx}\right)}{dx}+(G+\delta G)y_1 = 0,
\end{equation}
\ifTrans
Multiply the equations \eqref{eq:31} and \eqref{eq:32} by $y_1$ and $y$,
and subtract, we
\fi
\ifTrans
obtain
\fi
$$
y_1\tfrac{d}{dx}\left(L\tfrac{dy}{dx}\right)-y\tfrac{d}{dx}\left(L\tfrac{dy_1}{dx}\right)
-\delta Gyy_1 = 0.
$$
\ifTrans
Multiply by $dx$, and integrate between the limits $x_0$ and $X$; 
we
will have
\fi
$$
\int_{x_0}^X y_1\tfrac{d}{dx}\left(L\tfrac{dy}{dx}\right)dx
-\int_{x_0}^X y\tfrac{d}{dx}\left(L\tfrac{dy_1}{dx}\right)dx
-\int_{x_0}^X \delta Gyy_1dx = 0;
$$
\ifTrans
apply integration by parts to the first two terms, and 
if
$\delta h$
is infinitely small, the increase of $y$ will be too, and we will have, by replacing $y_1$, by $y+\delta y$,
\fi
\begin{equation}
\tag{A}\label{eq:3A}
\left\{
\begin{array}{l}
\left[L\left(\tfrac{dy}{dx}\delta y-y\delta\tfrac{dy}{dx}\right)\right]_X-
\left[L\left(\tfrac{dy}{dx}\delta y-y\delta\tfrac{dy}{dx}\right)\right]_{x_0}\\
\qquad\qquad\qquad-\displaystyle\int_{x_0}^X y^2\delta G dx = 0
\end{array}
\right.
\end{equation}

\ifTrans
Sturm then supposes that the quantity $L\tfrac{dx}{dy}: y$ undergoes for
$x=x_0$, by the increase of $h$, a variation of a given sign;
for the research we are proposing, imagine instead that $y$ is zero or
maximum or minimum for $x=x_0$, and whatever the value of the parameter $h$.
So for $x=x_0$, you have to make either $y=0,\delta y=0 $ or $\tfrac{dy}{dx}=0$, and $\delta\frac{dy}{dx}=0$; in both cases, the second bracket is zero, and it remains
\fi
\begin{equation}
\tag{3}\label{eq:33}
\left[L\left(\tfrac{dy}{dx}\delta y-y\delta\tfrac{dy}{dx}\right)\right]_X = 
\int_{x_0}^X y^2\delta G dx;
\end{equation}
\ifTrans
the first member is therefore of the same sign as the increase that
$G$ takes by the variation of the parameter.

We can now recognize if the roots in $x$ of the equation $y_1=0$ are
smaller or larger than those in the equation $y=0$. $y$ is a function 
of $x$ and $h$, and if for $x=X$ we have
\fi
$$ y=0, $$
\ifTrans
as soon as we give $h$ an increase $\delta h$, $y$ will no longer be
zero, unless we also give $x$ an increase $\delta x$ such that we have
\fi
$$
\tfrac{dy}{dx}\delta x + \tfrac{dy}{dh}\delta h = 0,
$$
\ifTrans
or
\fi
$$
\tfrac{dy}{dx}\delta x + \delta y = 0,
$$
\ifTrans
and therefore the root undergoes an increase equal to
\fi
\begin{equation}
\tag{4}\label{eq:34}
\delta x = -\tfrac{\delta y}{\tfrac{dy}{dx}}.
\end{equation}
\ifTrans
Suppose $L$ is positive; $y$ being zero for $x = X$, it follows from the
formula \eqref{eq:33} that if $\delta G$ is positive, 
$\tfrac{dy}{dx}\delta y$ is also positive, and, as a result of the
formula \eqref{eq:34}, that $\delta x$ is negative; so the roots of 
$y_1 = 0$ are smaller than those of $y = 0$, and similarly we see that 
if $\delta G$ is negative, the roots of $y_1 = 0$ are larger than those
of $y = 0$.

If we then imagine that we give to the parameter $h$, no longer an
infinitely small increase, but a finite increase, and that $h$ increases
from $h$ to $h_1$, if at the same time $G$ increases throughout this
interval, or decreasing throughout, the previous conclusions, relating to
the roots of $y = 0$ and $y_1 = 0$, are applicable, as can be recognized
by dividing the interval from $h$ to $h_1$ into infinitely small parts.

These considerations are due to Sturm (Journal of M. Liouville,
$1^{\rm {st}}$ series, t. I, p. 106); but we will show how they can be
used to recognize if a function is periodic, and we will get new results.

{\bf 7}. Let us first try to assume that $N$ is the square of an
integer $g^2$, and, substituting the letter $h$ for $\lambda c$, $P$ is
given by the equation
\fi
\begin{equation}
\tag{5}\label{eq:35}
\tfrac{d^2P}{d\alpha^2}+\left(g^2-4h^2\cos^2\alpha\right)P = 0.
\end{equation}
\ifTrans
It is very easy to recognize that the general solution of a linear
differential equation of the second order, such as the preceding ones,
can ordinarily be divided into two particular solutions, one of which is
zero, and the other of which is maximum or minimum for the value zero given to the variable. So let
\fi
$$
P = P_1+P_2,
$$
\ifTrans
$P_1$ being a solution which becomes zero [\emph{annule}] for $\alpha = 0$ and $P_2$ a solution which is maximum or minimum for this value.

These are the two functions $P_1$ and $P_2$ that we will examine. In the
case where $h$ becomes zero, they satisfy the equation
\fi
\begin{equation}
\tag{6}\label{eq:36}
\tfrac{d^2P'}{d\alpha^2}+g^2P' = 0,
\end{equation}
\ifTrans
and $P_1$ is reduced to $A\sin g\alpha$, $P_2$ to $B\cos g\alpha$.

$P_1$ becomes zero for $\alpha=0$, like $P'=A\sin g\alpha$; however
the coefficient of $P'$, in equation \eqref{eq:36}, is always greater
than that of $P$ in equation \eqref{eq:35}; it therefore follows from
what we have seen above that the roots of $P'= 0$ are smaller than those
of $P=0$; the roots of $P'=0 $ are from $0$ to $2\pi$,
\fi
$$
0,\quad \tfrac{\pi}{g},\quad \tfrac{2\pi}{g},\ldots,\quad \tfrac{(2g-1)\pi}{g}.
$$

\ifTrans
Let us give to $h$ an excessively small value, and divide, from the $x$-axis, a circumference which has its center at the origin into $2g$ equal parts, then lead to the points of division the rays $Oa$, $Ob$, $Oc$,\ldots; the roots of $P'=0$ are equal to the angles $aOb$, $aOc$, \ldots, and the roots of $P_1=0 $ are larger and represented by the angles $aOb'$, $aOc'$,\ldots. But when after a turn on the circumference we return to the point $a$, $P'$ again becomes zero at the value $\alpha=2\pi$, while $P_1$, which is zero for $\alpha=0 $, is not for $\alpha=2\pi$, but for a slightly larger value. $P_1$ does not therefore take the same value when we increase the arc $\alpha$ by a circumference.

To demonstrate that $P_2$, taken from equation \eqref{eq:35}, does not have $2\pi$ for period, apply the formula \eqref{eq:33}, replacing $x$ by $\alpha$, $L$ and $G$
\fi
\ifTrans
by $1$ and $g^2-4h^2\cos^2\alpha$, $y$ by $P_2$, $x_0$ and $X$ by $0$ and $2\pi$, and we will have
\fi
$$
\left(\tfrac{dP_2}{d\alpha}\delta P_2 -
P_2\delta\tfrac{dP_2}{d\alpha}\right)_{2\pi}
+4\left(2h\delta h+\delta h^2\right)\int_0^{2\pi}P_2^2\cos^2\alpha\, d\alpha=0,
$$
\ifTrans
a formula where we do not have to take account of $\delta h^2$ when $h$
is zero.
Suppose we vary $h$ from zero to the infinitely small value $\delta h$; 
$P_2$, for $h=0$, reduces to $B\cos g\alpha$, and it is then maximum for
$\alpha=2\pi$ as for $\alpha=0 $; so, by making $h=0$ in this formula, $\tfrac{dP_2}{d\alpha}$ becomes zero, and we have
\fi
$$
\left(P_2\delta\tfrac{dP_2}{d\alpha}\right)_{2\pi} = 
4\left(\delta h\right)^2\int_0^{2\pi}P_2^2\cos^2\alpha\, d\alpha;
$$
\ifTrans
the second member is essentially positive, so $\delta\tfrac{dP_2}{d\alpha}$ is not zero for $\alpha=2\pi$, or $\tfrac{dP_2}{d\alpha}$ is not zero for $\alpha=2\pi$ when we make $h=\delta h$; so finally $P_2$, which is maximum for $\alpha=0$, is not for $\alpha=2\pi $, and the function is not periodic.

If we take for the constant
\fi
$$
N=g^24 + h^2.
$$
\ifTrans
still designating by $g$ an integer, the equation which gives $P$ would become
\fi
$$
\tfrac{d^2P}{d\alpha^2}+\left(g^2+4h^2\sin^2\alpha\right)P=0.
$$

\ifTrans
By defining the particular solutions $P_1$ and $P_2$ as above, we will
recognize that the roots of $P_1=0$ and of $P_2=0$ are smaller than those
of $A\sin g\alpha=0$ and of $B\cos g\alpha=0$, and we will also
demonstrate, as before, that $P_1$ and $P_2$ are not periodic functions.

{\bf 8}. Finally, let us take for $N$ the expression
\fi
$$
N = g^2+2h^2,
$$
\ifTrans
we have the equation
\fi	
$$
\tfrac{d^2P}{d\alpha^2}+\left(g^2-2h^2\cos 2\alpha\right)P=0,
$$
\ifTrans
and we will demonstrate that $P$ is then a periodic function if $h$ is
excessively small, that is to say if we can neglect the powers of $h^2$
greater than the first.[\footnote{Nevertheless, if $g=1$, take $N=1+h^2$
or $=1+3h^2$, depending on whether it is $P_1$ or $P_2$.}]

If we consider only the solutions which are zero, or maximum or minimum
for $\alpha=0$, and which we have designated by $P_1$ and $P_2$, we have,
according to \eqref{eq:33}, the equation
\fi
\begin{equation}
\tag{7}\label{eq:37}
\tfrac{dP}{d\alpha}\delta P-P\delta\tfrac{dP}{d\alpha} =
-2\left(2h\delta h+\delta h^2\right)\int_0^{\alpha}P^2\cos 2\alpha\, d\alpha.
\end{equation}

\ifTrans
Instead of allowing just 
any $h$, let us take it equal to zero, and give
it the increase $\delta h$; then apply formula \eqref{eq:37} by making
$\alpha=2\pi$. If it is $P_1$ that we are considering, we will have
\fi
$$
P_1^2\cos 2\alpha=\sin^2g\alpha\cos2\alpha=
\tfrac{\cos2\alpha}{2}-
\tfrac{\cos2(g+1)\alpha+\cos2(g-1)\alpha}{4},
$$
\ifTrans
and if $g$ is not equal to $1$, we will have
\fi
$$
\int_0^{\tfrac{\pi}{2}}P_1^2\cos2\alpha\, d\alpha=0.
$$

\ifTrans
If we consider $P_2$, we have
\fi
$$
P_2^2\cos 2\alpha=\cos^2g\alpha\cos2\alpha=
\tfrac{\cos2\alpha}{2}-
\tfrac{\cos2(g+1)\alpha+\cos2(g-1)\alpha}{4},
$$
\ifTrans
and if g is not equal to $1$, we still have
\fi
$$
\int_0^{\tfrac{\pi}{2}}P_2^2\cos2\alpha\, d\alpha=0;
$$
\ifTrans
so equation \eqref{eq:37} is reduced in both cases to
\fi
\begin{equation}
\tag{8}\label{eq:38}
\left(\tfrac{dP}{d\alpha}\delta P-P\delta\tfrac{dP}{d\alpha}\right)_{\tfrac{\pi}{2}}=0;
\end{equation}
\ifTrans
but for $h=0$, $P_1$ and $P_2$ become $\sin g\alpha$ and $\cos g\alpha$,
and one of the two functions is zero and the other maximum for 
$\alpha=\tfrac{\pi}{2}$; now from this formula we conclude that if $P$
is zero for $\alpha=\tfrac{\pi}{2}$, $\delta P$ is too, and that if 
$\tfrac{dP}{d\alpha}$ is zero for $\alpha=\tfrac{\pi}{2}$, 
$\delta\tfrac{dP}{d\alpha}$ is zero at the same time. So, for a very
small value of $h$ and for $\alpha=\tfrac{\pi}{2}$, $P_1$ is zero or
maximum like $\sin g\alpha$, and $P_2$ is zero or maximum like 
$\cos g\alpha$.

Suppose now that $h$ is no longer excessively small, but has some value;
and, setting $N=R+2h^2$, consider the equation
\fi
\begin{equation}
\tag{9}\label{eq:39}
\tfrac{d^2P}{d\alpha^2}+\left(R-2h^2\cos2\alpha\right)P=0,
\end{equation}
\ifTrans
in which $R$ depends on $h$, and is reduced to the square $g^2$ of an
integer for $h=0$; by applying the equation \eqref{eq:33}, we have
\fi
$$
\tfrac{dP}{d\alpha}\delta P-P\delta\tfrac{dP}{d\alpha} = 
\int_0^{\alpha}P^2\left(\delta R-4h\delta h\cos2\alpha\right)d\alpha;
$$
\ifTrans
so imagine that we know how to determine the constant $R$ so that the
integral
\fi
$$
\int_0^{\tfrac{\pi}{2}}P^2\left(\delta R-4h\delta h\cos2\alpha\right)d\alpha
$$
\ifTrans
is zero, regardless of $h$: the property we just obtained when $h$ is
very small takes place for any value of $h$, because we will still have
the equation \eqref{eq:38}, and $P$ will be zero or maximum for 
$\alpha=\tfrac{\pi}{2}$, depending on whether $\sin g\alpha$ and 
$\cos g\alpha$, to which it reduces for $h=0$, is zero or maximum.
\fi

\ifTrans
Note now that nothing indicates that, for the same value of $g$, the
constant $R$ is the same in the functions $P_1$, and $P_2$; it is indeed
different, and the general solution of equation \eqref{eq:39} cannot be
periodic. To fix the ideas, let us choose the constant so that we will have
\fi
$$
\int_0^{\tfrac{\pi}{2}}P_1^2\left(\delta R-4h\cos2\alpha\right)d\alpha=0,
$$
\ifTrans
and I say that the function $P_1$ 
will take the same values or equal
values and of opposite sign in each quadrant, so that $\alpha$ being
between $0$ and $\tfrac{\pi}{2}$ the four quantities
\fi
$$
P_1(\alpha),\quad P_1(\pi-\alpha),\quad P_1(\pi+\alpha),\quad P_1(2\pi-\alpha)
$$
\ifTrans
are equal, except for the sign, and that $P_1$ is periodic.

As we will have occasion to see later, if we set $\nu=\cos\alpha$, the
general solution of the differential equation which gives $P$ is the sum
of two particular solutions which develop thus
\fi
$$
P' = A_0+A_1\nu^2+A_2\nu^4+A_3\nu^6+\ldots,
$$
$$
P'' = B\nu+B_1\nu^3+B_2\nu^5+B_3\nu^7+\ldots.
$$

\ifTrans
The first is maximum and the second zero for $\nu=0$, 
wherefore 
$\alpha=\tfrac{\pi}{2}$.

It follows from this that $P_1$ is equal to $P'$ or $P''$, depending on
whether it is zero or maximum for $\alpha=\tfrac{\pi}{2}$; however, if we
change $\nu$ to $-\nu$ or $\alpha$ to $\pi-\alpha$, $P'$ remains
constant and $P''$ changes sign only; therefore $P$, remains the same,
except for the sign, when we replace $\alpha$ by $\pi-\alpha$.

To go to the third quadrant, we note that the general solution of $P$ can
be divided into two solutions, one of which is even, and the other of
which is odd, according to the powers of $\nu'=\sin\alpha$; $P_1$, which
is zero for $\alpha=\pi$, is equal to the odd solution in $\nu'$, and we
conclude
\fi
$$
P_1(\pi+\alpha)=-P_1(\pi-\alpha),\quad P_1(\pi+\alpha)=\pm P_1(\alpha).
$$
\ifTrans
Finally, the values of $P_1$ in the fourth quadrant can be obtained 
in the same way. So the function $P_1$ takes the same values down to the
sign in each quadrant, and behaves in sign changes like $\sin g\alpha$,
and it has $2\pi$ for period.

If we determine the constant $R$ so that we have
\fi
$$
\int_0^{\frac{\pi}{2}}P_2^2\left(\tfrac{dR}{dh}-4h\cos2\alpha\right)d\alpha=0,
$$
\ifTrans
we arrive at similar conclusions for $P_2$, which behaves in the passage
from one quadrant to the next like $\cos g\alpha$. It is therefore only
necessary to study the functions $P_1$ and $P_2$ between the limits 
$\alpha=0$ and $\alpha=\tfrac{\pi}{2}$.

{\bf 9}. If we first regard 
$h$ as very small, the constant $R$ is
reduced to almost $g^2$, and the differential equation to
\fi
$$
\tfrac{d^2P}{d\alpha^2}+\left(g^2-2h^2\cos2\alpha\right)P=0.
$$

\ifTrans
If we give $h$ the increase $\delta h$, any root of
$$
P_1=0\quad\rm{or\ of}\quad P_2=0
$$
undergoes a variation whose value is
\fi
$$
\delta\alpha=-\delta P:\tfrac{dP}{d\alpha}.
$$

\ifTrans
We have the general formula
\fi
\begin{equation}
\tag{10}\label{eq:310}
\left\{
\begin{array}{l}
\tfrac{dP}{d\alpha}\delta P-P\delta\tfrac{dP}{d\alpha}\\
\quad =\left(\tfrac{dP}{d\alpha}\delta P-P\delta\tfrac{dP}{d\alpha}\right)_a
-2\left(2h\delta h+\delta h^2\right)\displaystyle\int_a^{\alpha}P^2\cos2\alpha\, d\alpha.
\end{array}
\right.
\end{equation}
\ifTrans
Let us make $a=0$; suppose that $P$ represents $P_1$ or $P_2$, and that
$\alpha$ is a root of $P=0$, enclosed between $0$ and $\tfrac{\pi}{4}$,
the preceding equation
\fi
\ifTrans
becomes
\fi
$$
\tfrac{dP}{d\alpha}\delta P =
-2\left(2h\delta h+\delta h^2\right)\int_0^{\alpha}P^2\cos2\alpha\, d\alpha.
$$

\ifTrans
All the elements of the integral are positive; therefore 
$\tfrac{dP}{d\alpha}\delta P$ is negative and the variation of the roots
positive, when we give to $h$ the increase $\delta h$.

Let us make $a=\tfrac{\pi}{2}$, and suppose that $\alpha$ is now a root
of $P=0$, enclosed between $\tfrac{\pi}{4}$ and $\tfrac{\pi}{2}$, we
deduce from the same formula
\fi
$$
\tfrac{dP}{d\alpha}\delta P =
2\left(2h\delta h+\delta h^2\right)\int_0^{\tfrac{\pi}{2}}P^2\cos2\alpha\, d\alpha.
$$
\ifTrans
that the second member is negative; so the root increase is still
positive.

Suppose for example that it is $P_2$, and that $g$ is
even; then the function $P$ is maximum like $\cos g\alpha$ for 
$\alpha=0$ and $\alpha=\tfrac{\pi}{2}$. If $Ob$, $Oc$, $Od$, \ldots make
angles with $Ox$ equal to the roots of the equation $\cos g\alpha=0$,
these lines can represent nodal lines of the circular membrane, and the
arcs $bc$, $cd$,\ldots are equal to each other and double the extreme arcs 
$ab$ and $ef$ of the quadrant $af$.

Let us consider an ellipse whose eccentricity is very small, and take
[\emph{menons}] the asymptotes of the hyperbolic nodal lines 
$Ob'$, $Oc'$, $Od'$,\ldots; it follows from what we have shown that the
angles $aOb'$, $aOc'$, $aOd'$,\ldots are respectively larger than $aOb$,
$aOc$,\ldots. But there is more: the angles $b'Oc'$, $c'Od'$,\ldots are
$>bOc$, $cOd$,\ldots in the first half of the quadrant and are less in
the second half.

To prove it, let us denote by $\alpha_1$ and $\alpha_2$, two consecutive
roots of the equation $P_2=0$, and consider the function
\fi
$$
\Pi = A\sin g(\alpha-\alpha_1),
$$
\ifTrans
which becomes zero for $\alpha=\alpha_1$; $P_2$ is not reduced to $\Pi$
for $h=0$; but we can imagine a function $P$ which satisfies the
differential equation
\fi
\ifTrans
of the second order, and which by the variation of $h$ passes from $\Pi$
to $P_1$ while remaining constantly zero for $\alpha=\alpha_1$. So for 
$\alpha=\alpha_1$, we will have
\fi
$$
P=0,\quad \delta P=0,
$$
\ifTrans
and by making $\alpha=\alpha_2$, and $a=\alpha_1$ in equation 
\eqref{eq:310}, we get
\fi
$$
\left(\tfrac{dP}{d\alpha}\delta P\right)_{\alpha=\alpha_1} = 
-2\left(2h\delta h+\delta h^2\right)
\int_{\alpha_1}^{\alpha_2}P^2\cos2\alpha\, d\alpha;
$$
\ifTrans
if $\alpha_1$ and $\alpha_2$ are less than $\tfrac{\pi}{4}$, the integral
of the second member is positive, and the first member is negative; so
the variation of the root $\alpha$,
\fi
$$
\delta\alpha_2=-\delta P:\tfrac{dP}{d\alpha},
$$
\ifTrans
but this time counted from $\alpha_1$, is positive; therefore the
interval of the two roots $\alpha_1$ and $\alpha_2$ has increased; it is
therefore larger than $\tfrac{2\pi}{g}$. We would see that on the
contrary if $\alpha_1$ and $\alpha_2$ are between $\tfrac{\pi}{4}$ and 
$\tfrac{\pi}{2}$, the integral is negative, and that the interval between
two roots decreases when $h$ increases, while keeping very small values.

{\bf 10}. But whatever the size of $h$ and whatever the direction in
which the constant $R$ varies, when we give an infinitely small increase
to $h$, the roots undergo infinitely small modifications, and those which
were included between $0$ and $\tfrac{\pi}{2}$ and 
will remain constant;
because $Oa$ and $Of$ are lines where $P_2$ is maximum and cannot become
zero, and that consequently these roots by changing magnitude cannot
cross. Then this property belongs in all cases to the functions $P_1$ and
$P_2$, the first of which is zero, and the second maximum for $\alpha=0$,
while one becomes zero and the other is maximum for 
$\alpha=\tfrac{\pi}{2}$, according to the parity of $g$. So, for example,
if $P_2$ becomes zero for $\alpha=\tfrac{\pi}{2}$, it is impossible that by
increasing $h$ a root  
\fi
\ifTrans
between $0$ and $\tfrac{\pi}{2}$ crosses the limit $\tfrac{\pi}{2}$;
because if for a value of $h$ a root between $0$ and $\tfrac{\pi}{2}$
becomes equal to $\tfrac{\pi}{2}$, the equation $P_2=0$ would have a
double root for $\alpha=\tfrac{\pi}{2}$; therefore $P$ and 
$\tfrac{dP}{d\alpha}$, and consequently the derivatives of all orders,
would become zero together for the same value of $\alpha $; which is
impossible.

Now for $h=0$ the functions $P_1$ and $P_2$ are reduced to 
$\sin g\alpha$ and $\cos g\alpha$, and become zero $g$ times from $0$ to
$\pi$; therefore whatever $h$, the equations $P_1=0$ and $P_2=0$ also
have $g$ roots from $0$ to $\pi$ (admitting among these roots that which
would be zero, but not that which would be equal to $\pi$).

\begin{center}
\emph{Development of the functions $P_1$ and $P_2$ according to 
powers of $h$.}
\end{center}

{\bf 11}. To develop according to powers of $h$ the solutions $P_1$ and $P_2$ of the equation
\fi
$$
\tfrac{d^2P}{d\alpha^2} + \left(R-2h^2\cos2\alpha\right)P = 0
$$
\ifTrans
that have $2\pi$ for period, and the first of which is zero, the second maximum for $\alpha=0$, we will let
\fi
$$
R = g^2 + \beta h^4 + \gamma h^6 + \delta h^8 +\ldots,
$$
\ifTrans
by designating by $g$ any integer, and we will seek to determine the coefficients according to the condition that $P_1$ and $P_2$ are periodic.

First consider $ P_2 $; let in the differential equation
\fi
$$
P=P_2=\cos g\alpha + h^2{\rm p},\quad R = g^2+Bh^4,
$$
\ifTrans
and we will have
\fi
$$
\tfrac{d^2{\rm p}}{d\alpha^2} + (g^2 - 2h^2\cos2\alpha + Bh^4){\rm p} -
\left(2\cos2\alpha\cos g\alpha - Bh^2\cos g\alpha\right) = 0.
$$

\ifTrans
Then let
\fi
$$
{\rm p}=p+h^2{\rm p}_1,
$$
\ifTrans
and we will have
\fi
\begin{equation}
\tag{I}\label{eq:4I}
0=\tfrac{d^2p}{d\alpha^2}+g^2p-2\cos2\alpha\cos g\alpha,
\end{equation}
\begin{equation}
\tag{$b$}\label{eq:4b}
\left\{
\begin{array}{p{0.1em}p{0.2em}p{15em}}
$0$&$=$&$\tfrac{d^2{\rm p}_1}{d\alpha^2}+
\left(g^2-2h^2\cos2\alpha+Bh^4\right){\rm p}_1$\\
&&$+\left(-2\cos2\alpha + Bh^2\right)p+B\cos g\alpha.$
\end{array}
\right.
\end{equation}

\ifTrans
To solve equation \eqref{eq:4I}, we will replace 
$2\cos2\alpha\cos g \alpha$ by $\cos(g+2)\alpha+\cos(g-2)\alpha$, and we will let
\fi
$$
p=a\cos(g + 2)\alpha + b\cos g\alpha +c\cos(g-2)\alpha;
$$
\ifTrans
we find immediately
\fi
$$
a=\tfrac{-1}{4(g+1)},\quad c=\tfrac{1}{4(g-1)};
$$
\ifTrans
for $b$, it is not determined, and indeed $P_2$ is reduced to 
$\cos g\alpha$ for $h=0$; but if we suppose that we have obtained its
expression, and that we multiply it by $1+Bh^2$, this new expression can
still represent $P_2$, and the coefficient $b$ thereby changes in an
arbitrary way.

Since we can give $b$ the value we want, we will make
\fi
$$
b = 0.
$$

\ifTrans
In equation \eqref{eq:4b}, we let
\fi
$$
{\rm p}_1 = p_1 + h^2{\rm p}_2,\quad B=\beta + Ch^2,
$$
\ifTrans
and we will have the two equations
\fi
\begin{equation}
\tag{II}\label{eq:4II}
\tfrac{d^2p_1}{d\alpha^2}+g^2p_1-2\cos2\alpha p+\beta\cos g\alpha = 0,
\end{equation}
\begin{equation}
\tag{$c$}\label{eq:4c}
\left\{
\begin{array}{p{0.1em}p{0.2em}p{18em}}
$0$&$=$&$\tfrac{d^2{\rm p}_2}{d\alpha^2}+
\left(g^2-2h^2\cos2\alpha+\beta h^4 + Ch^6\right){\rm p}_2$\\
&&$+\left(-2\cos2\alpha + \beta h^2 + Ch^4\right)p_1$\\
&&$+\left(\beta + Ch^2\right)p + C\cos g\alpha.$
\end{array}
\right.
\end{equation}
\ifTrans
To solve equation \eqref{eq:4II}, we must replace $2\cos2\alpha p$ by its
\fi
 value
\begin{multline}
\qquad\qquad a\cos(g + 4)\alpha + b\cos(g +2)\alpha\\
+(a+c)\cos g\alpha + b\cos(g-2)\alpha + c\cos(g-4)\alpha,\qquad\qquad
\end{multline}
\ifTrans
then substitute for $p_1$
\fi
\begin{multline}
\qquad\qquad p_1 = d\cos(g + 4)\alpha + e\cos(g + 2)\alpha\\
+ f\cos g\alpha + h\cos(g-2)\alpha + k\cos(g-4)\alpha,\qquad\qquad
\end{multline}
\ifTrans
and we find, by equating to zero the coefficients of the cosines of the
different arcs,
\fi
$$
d = -\tfrac{a}{8(g+2)},\quad e=\tfrac{b}{4(g+1)},\quad 
h=\tfrac{b}{4(g-1)},\quad k=\tfrac{c}{8(g-2)},
$$
$$
\beta=a+c.
$$
\ifTrans
$b$ is left arbitrary in these formulas; if we assume $b=0$, we have
\fi
$$
d=\tfrac{1}{32(g+1)(g+2)},\quad e=0,\quad h=0,\quad k=\tfrac{1}{32(g-1)(g-2)},
$$
$$
\beta=\tfrac{1}{2(g^2-1)}.
$$

\ifTrans
The coefficient $ f $ remains undetermined, and indeed, if we imagine
that we have obtained the expression of $P_2$ and that we multiply it by
$1+Bh^2+Ch^4$, we will not only change the coefficient $b$ in an
arbitrary manner, but also the coefficient $f$; the simplest thing is to
make $f=0$.

In equation \eqref{eq:4c}, let
\fi
$$
{\rm p}_2 = p_2+h^2{\rm p}_3,\quad C=\gamma+Dh^2,
$$
\ifTrans
and we have
\fi
\begin{equation}
\tag{III}\label{eq:4III}
\tfrac{d^2p_2}{d\alpha^2}+g^2p_2-2\cos2\alpha p_1+\gamma\cos g\alpha = 0,
\end{equation}
\begin{equation}
\tag{$d$}\label{eq:4d}
\left\{
\begin{array}{p{0.1em}p{0.2em}p{18em}}
$0$&$=$&$\tfrac{d^2{\rm p}_3}{d\alpha^2}+
\left(g^2-2h^2\cos2\alpha+\beta h^4 + Ch^6\right){\rm p}_3$\\
&&$+\left(-2\cos2\alpha + \beta h^2 + Ch^4\right)p_2$\\
&&$+\left(\beta + Ch^2\right)p_1 + Cp + D\cos g\alpha.$
\end{array}
\right.
\end{equation}
\ifTrans
Let us replace in equation \eqref{eq:4III} $2\cos2\alpha p_1$ by
\fi
\begin{multline}
d\cos(g + 6)\alpha + e\cos(g + 4)\alpha +(d+f)\cos(g + 2)\alpha
+(e + h)\cos g\alpha\\
+(k+f)\cos(g-2)\alpha + h\cos(g-4)\alpha + k\cos(g-6)\alpha,
\end{multline}
\ifTrans
and let
\fi
\begin{multline}
\quad p_2 = l\cos (g + 6)\alpha + m\cos(g + 4)\alpha + n\cos(g + 2)\alpha\\
+ \varpi\cos g\alpha + q\cos(g-2)\alpha + r\cos(g-4)\alpha + 
s\cos (g-6)\alpha;\quad 
\end{multline}
\ifTrans
we will have
\fi
$$
l=\tfrac{-d}{12(g+3)},\quad m=\tfrac{-e}{8(g+2)},\quad 
r=\tfrac{h}{8(g-2)},\quad s=\tfrac{k}{12(g-3)},
$$
$$
n=\tfrac{-d-f-\beta\alpha}{4(g+1)},\quad q=\tfrac{k+f-\beta c}{4(g-1)},
\quad \gamma=e+h-\beta b = 0.
$$

\ifTrans
These are the expressions of $l$, $m$, $n$,\ldots, regardless of the
values given to $b$ and $f$; and if we suppose them zero, we get
\fi
$$
l=\tfrac{-1}{2^7\cdot 3(g+1)(g+2)(g+3)},\quad m=0,\quad r=0,\quad s=\tfrac{1}{2^7\cdot 3(g-1)(g-1)(g-3)},
$$
$$
n=\tfrac{-(g^2+4g+7)}{2^7(g+1)^3(g-1)(g+2)},\quad 
q=\tfrac{g^2-4g+7}{2^7(g-1)^3(g+1)(g-2)},\quad \gamma=0.
$$
\ifTrans
$\varpi$ is undetermined, like $b$ and $f$, and we will make it zero too.

Now that we see what kind of simplification 
results from the
hypothesis of the nullity of arbitrary constants, and that we recognize
that it brings about the disappearance of the terms of even rank in $p$,
$p_1$, $p_2$, etc., immediately make these reductions in the calculations
that follow. Let us put in equation \eqref{eq:4d}
\fi
$$
{\rm p}_3 = p_3 + h^2{\rm p}_4,\quad D=\delta + E h^2,
$$
\ifTrans
we will get the two equations
\fi
\begin{equation}
\tag{IV}\label{eq:4IV}
\tfrac{d^2p_3}{d\alpha^2}+g^2p_3-2\cos2\alpha\, p_2+\beta p_1+\delta\cos g\alpha = 0,
\end{equation}
\begin{equation}
\tag{$e$}\label{eq:4e}
\left\{
\begin{array}{p{0.1em}p{0.2em}p{22em}}
$0$&$=$&$\tfrac{d^2{\rm p}_4}{d\alpha^2}+
\left(g^2-2h^2\cos2\alpha+\beta h^4 + Ch^6\right){\rm p}_4$\\
&&$+\left(-2\cos2\alpha + \beta h^2 + Ch^4\right)p_3+(\beta+Ch^2)p_2$\\
&&$+Ch^2p_1 + \left(\delta + Eh^2\right)p + E\cos g\alpha.$
\end{array}
\right.
\end{equation}
\ifTrans
Let us replace in equation \eqref{eq:4IV} $2\cos\alpha p_2$ by
\fi
\begin{multline}
\qquad\quad l\cos(g+8)\alpha+(l+n)\cos(g+4)\alpha+(q+n)\cos g\alpha\\
+(q+s)\cos(g-4)\alpha+s\cos(g-8)\alpha,\qquad\quad 
\end{multline}
\ifTrans
and let
\fi
\begin{multline}
\qquad\quad p_3=R_1\cos(g+8)\alpha+R_2\cos(g+4)\alpha+R_3\cos g\alpha\\
+R_4\cos(g-4)\alpha+R_5\cos(g-8)\alpha;\qquad\quad 
\end{multline}
\ifTrans
we will have
\fi
$$
R_1=\tfrac{-l}{16(g+4)},\quad R_2=\tfrac{-(l+n)+\beta d}{8(g+2)},\quad
R_4=\tfrac{q+s-\beta k}{8(g-2)},\quad R_5=\tfrac{s}{16(g-4)},
$$
$$
\delta=q+n,
$$
\ifTrans
or, by performing the calculations,
\fi
$$
R_1=\tfrac{1}{2^{11}\cdot 3(g+1)(g+2)(g+3)(g+4)},\quad 
R_5=\tfrac{1}{2^{11}\cdot 3(g-1)(g-2)(g-3)(g-4)},
$$
$$
R_2=\tfrac{g^3+7g^2+20g+20}{2^8\cdot 3(g+1)^3(g-1)(g+2)^2(g+3)},\quad 
R_5=\tfrac{g^3-7g^2+20g-20}{2^8\cdot 3(g-1)^3(g+1)(g-2)^2(g-3)},
$$
$$
\delta=\tfrac{5g^2+7}{32(g^2-1)^3(g^2-2^2)};\quad
$$
\ifTrans
add to these values $R_3=0$.

If we again let
\fi
$$
{\rm p}_4 = p_4 + h^2{\rm p}_5,\quad E=\varepsilon+Hh^2,
$$
\ifTrans
$p_4$ will be given to us by the equation
\fi
$$
\tfrac{d^2 p_4}{d\alpha^2} + g^2p_4 - 2\cos2\alpha\, p_3 + \beta p_2\\
+\delta p + \varepsilon\cos g\alpha = 0,
$$
\ifTrans
and we will have
\fi
\begin{multline}
\qquad p_4=S_1\cos(g+10)\alpha+S_2\cos(g+6)\alpha+S_3\cos(g+2)\alpha\\
-S_4\cos(g-2)\alpha+S_5\cos(g-6)\alpha+S_6\cos(g-2)\alpha,\qquad
\end{multline}
\ifTrans
taking
\fi
$$
S_1=\tfrac{-R_1}{4\cdot 5(g+5)},\quad
S_2=\tfrac{-R_1-R_2+\beta l}{4\cdot 3(g+3)},\quad
S_3=\tfrac{-R_2+\beta n+\delta a}{4(g+1)},
$$
$$
S_4=\tfrac{R_4-\beta q-\delta c}{4(g-1)}\quad
S_5=\tfrac{R_4+R_5-\beta s}{4\cdot 3(g-3)},\quad
S_6=\tfrac{R_5}{4\cdot 5(g-5)},
$$
$$
\varepsilon=0
$$
\ifTrans
or
\fi
$$
S_1=\tfrac{-1}{2^{11}\cdot3\cdot4\cdot5(g+1)(g+2)(g+3)(g+4)(g+5)},\ \,
S_6=\tfrac{1}{2^{11}\cdot3\cdot4\cdot5(g-1)(g-2)(g-3)(g-4)(g-5)},
$$
$$
S_2=\tfrac{-(g^4+11g^3+49g^2+101g+78)}{2^{13}(g+1)^3(g+2)^2(g+3)^2(g+4)(g-1)},\ \,
S_2=\tfrac{g^4-11g^3+49g^2-101g+78}{2^{13}(g-1)^3(g-2)^2(g-3)^2(g-4)(g+1)},
$$
$$
S_3=\tfrac{-(g^7+7g^6+18g^5+24g^4+63g^3+81g^2+206g+464)}{2^{10}\cdot3(g+1)^5(g+2)^2(g+3)(g-1)^3(g-2)},
$$
$$
S_4=\tfrac{g^7-7g^6+18g^5-24g^4+63g^3-81g^2+206g-464}{2^{10}\cdot3(g-1)^5(g-2)^2(g-3)(g+1)^3(g+2)}.
$$ 

\ifTrans
To complete this calculation, note that the coefficient $\eta$ of $h^{12}$ in the constant has the value
\fi
$$
\eta=S_3+S_4,
$$
\ifTrans
and, replacing $S_3$ and $S_4$ with their expressions,
\fi
$$
\eta=\tfrac{9g^6+22g^4-203g^2-116}{2^6(g^2-1)^5(g^2-2^2)^2(g^2-3^2)}.
$$

\ifTrans
So, by putting in $R$ the values of the first terms, we get
\fi
\begin{equation}
\tag{$l$}\label{eq:4l}
\left\{
\begin{array}{p{0.1em}p{0.2em}p{18em}}
$R$&$=$&$g^2+\tfrac{1}{2(g^2-1)}h^4+\tfrac{5g^2+7}{32(g^2-1)^3(g^2-4)}h^8$\\
&&$+\tfrac{9g^6+22g^4-203g^2-116}{64(g^2-1)^5(g^2-4)^2(g^2-9)}h^{12}+\ldots.$
\end{array}
\right.
\end{equation}

\ifTrans
{\bf 12}. It is time to notice that we cannot continue the development of
$P_2$ and the constant $R$ in this way without worrying about the value
of the integer $g$; because the coefficient of $h^4$ contains in its
denominator the factor $g-1$, the coefficient of $h^8$ the factor $g-2$,
the coefficient of $ h^{12}$ the factor $g-3$, and so on; so that whatever
\fi
\ifTrans
the integer taken for $g$, we will end up finding an infinite term. We
must even stop the development of $R$ before meeting an infinite term;
because, for a term of the constant to be accepted, the term must itself be of the same order as $P_2$.

For clarity, consider a special case, that of $g=4$ for example. The
coefficients of $h^8$ and $h^{12}$ in $R$ keep a finite value 
but 
must
however be rejected. To recognize 
this, 
let us resume the calculation of 
$p_3$, which needs to be modified, because the value of $R_5$ which
appears there becomes infinite.

The expression of $2\cos2\alpha p_2$, becomes
\fi
$$
l\cos12\alpha + (l+n)\cos8\alpha + (q+n+s)\cos4\alpha+q+s,
$$
\ifTrans
and the terms in $\cos g\alpha$ and $\cos (g-8)\alpha$ come together in one, in $\cos4\alpha$.

We will substitute in the equation \eqref{eq:4IV}
\fi
$$
p_3=R_1\cos12\alpha+R_2\cos8\alpha+R_3\cos4\alpha+R_4,
$$
\ifTrans
to determine the coefficients; but in the result of the substitution, the
coefficient of $\cos4\alpha$ having to be zero, we get
\fi
$$
\delta=q+n+s;
$$
\ifTrans
the value of $\delta$ must therefore be increased by the quantity $s$; we
have
\fi
$$
q+n=\tfrac{87}{2^7\cdot3^4\cdot5^3},\quad s=\tfrac{1}{2^8\cdot 3^2},
$$
\ifTrans
and consequently
\fi
$$
\delta=\tfrac{433}{2^8\cdot 3^3\cdot 5^3}.
$$

\ifTrans
We can clearly see how we would continue this development.

If $g$ is $>4$, the first three terms of $R$ are those found in the
expression \eqref{eq:4l}; if $g$ is $>6$, we must still take in the
development \eqref{eq:4l} the term in $h^{12}$, and so on.

We will consider the development of $P_2$ and the first terms of $R$ when
$g$ is less than $4$.
\fi

\ifTrans
If $g$ is zero, the expansion we found for $P_2$ is applicable, and it is
even the only case where we can apply it as far as we want.

If $g=2$, we obtain, by a special calculation,
\fi
$$
\begin{array}{p{0.4em}p{0.2em}p{17em}}
$P_2$&$=$&$\cos2\alpha
+h^2\left(-\tfrac{1}{12}\cos4\alpha+\tfrac{1}{4}\right)
+\tfrac{h^4}{384}\cos6\alpha$\\
&&$-h^6\left(\tfrac{1}{23040}\cos8\alpha
+\tfrac{43}{13824}\cos4\alpha+\tfrac{5}{192}\right)$\\
&&$+h^8\left(\tfrac{1}{2211840}\cos10\alpha+\tfrac{287}{2211840}\cos6\alpha\right)$\\
&&$+h^{10}\left(\tfrac{-1}{309657600}\cos12\alpha
-\tfrac{41}{16588800}\cos8\alpha\right.$\\
&&$\qquad\quad+\left.\tfrac{21059}{79626240}\cos4\alpha
+\tfrac{1363}{221184}\right)+\ldots;$
\end{array}
$$
\begin{equation}
\tag{A}\label{eq:4A}
R=4+\tfrac{5}{12}h^4-\tfrac{763}{13824}h^8
+\tfrac{1002419}{79626240}h^{12}+\ldots
\end{equation}

\ifTrans
If $g=4$, we have
\fi
$$
\begin{array}{p{0.4em}p{0.2em}p{26em}}
$P_2$&$=$&$\cos4\alpha
+h^2\left(-\tfrac{1}{20}\cos6\alpha+\tfrac{1}{12}\cos2\alpha\right)
+\tfrac{h^4}{384}\left(\tfrac{1}{960}\cos8\alpha+\tfrac{1}{192}\right)$\\
&&$+h^6\left(\tfrac{-1}{80640}\cos10\alpha
-\tfrac{13}{96000}\cos6\alpha+\tfrac{11}{17280}\cos2\alpha\right)$\\
&&$+h^8\left(\tfrac{1}{10321920}\cos12\alpha
+\tfrac{23}{6048000}\cos8\alpha-\tfrac{1}{92160}\right)$\\
&&$-h^{10}\left(\tfrac{1}{1857945600}\cos14\alpha
+\tfrac{53}{1032192000}\cos10\alpha\right.$\\
&&$\qquad\quad+\left.\tfrac{4037}{2419200000}\cos6\alpha
+\tfrac{439}{62208000}\cos2\alpha\right)+\ldots;$
\end{array}
$$
\begin{equation}
\tag{B}\label{eq:4B}
R=16+\tfrac{1}{30}h^4+\tfrac{433}{864000}h^8
-\tfrac{189983}{21772800000}h^{12}+\ldots
\end{equation}
\ifTrans
The expansions \eqref{eq:4A} and \eqref{eq:4B} contain only even powers
of $h^4$, and we will demonstrate below that $R$ enjoys this property
whenever $g$ is even.
\fi
\ifTrans
For $g=1$, we have the following formulas:
\fi
$$
\begin{array}{p{0.4em}p{0.2em}p{22em}}
$P_2$&$=$&$\cos\alpha-\tfrac{h^2}{8}\cos3\alpha
+h^4\left(\tfrac{1}{192}\cos5\alpha-\tfrac{1}{64}\cos3\alpha\right)$\\
&&$-h^6\left(\tfrac{1}{9216}\cos7\alpha
-\tfrac{1}{1152}\cos5\alpha+\tfrac{1}{1536}\cos3\alpha\right)$\\
&&$+h^8\left(\tfrac{1}{737280}\cos9\alpha
-\tfrac{1}{49152}\cos7\alpha\right.$\\
&&$\qquad\quad+\left.\tfrac{1}{24576}\cos5\alpha
+\tfrac{11}{36864}\cos3\alpha\right)+\ldots;$
\end{array}
$$
$$
R=1+h^2-\tfrac{1}{8}h^4-\tfrac{1}{64}h^6-\tfrac{1}{1536}h^8+\tfrac{11}{36864}h^{10}+\ldots.
$$
\ifTrans
For $g=3$, we have
\fi
$$
\begin{array}{p{0.4em}p{0.2em}p{20em}}
$P_2$&$=$&$\cos3\alpha
+h^2\left(-\tfrac{1}{16}\cos5\alpha+\tfrac{1}{12}\cos\alpha\right)$\\
&&$+h^4\left(\tfrac{1}{640}\cos7\alpha+\tfrac{1}{64}\cos\alpha\right)$\\
&&$+h^6\left(\tfrac{-1}{46080}\cos9\alpha
-\tfrac{7}{20480}\cos5\alpha
+\tfrac{1}{768}\cos\alpha\right)$\\
&&$+h^8\left(\tfrac{1}{2^{14}\cdot 3^2\cdot 5\cdot7}\cos11\alpha
-\tfrac{17}{2^{15}\cdot3^2\cdot5}\cos7\alpha\right.$\\
&&$\qquad\quad-\left.\tfrac{1}{2^{14}}\cos5\alpha
-\tfrac{1}{2^{13}}\cos\alpha\right)+\ldots;$
\end{array}
$$
$$
R=9+\tfrac{1}{16}h^4+\tfrac{1}{64}h^6+\tfrac{59}{61440}h^8
-\tfrac{3}{16384}h^{10}+\ldots.
$$
\ifTrans
{\bf 15}. We now propose to develop $P_1$, and as for the same value of
$g$ the constant $R$ has a different value in $P_1$ and $P_2$, let us now
represent it by $R'$. So we have the differential equation
\fi
\begin{equation}
\tag{$m$}\label{eq:4m}
\tfrac{d^2P_1}{d\alpha^2}+\left(R'-2h^2\cos2\alpha\right)P_1=0,
\end{equation}
\ifTrans
and you have to find the solution that becomes zero for $\alpha=0$ and choose
\fi
$$
R'=g^2+\beta h^4+\gamma h^6+\delta h^8+\varepsilon h^{10}+\eta h^{12}+\ldots
$$
\ifTrans
so that it is periodic. Let
\fi
$$
P_1=\sin g\alpha+h^2p+h^4p_1+h^6p_2+h^8p_3+\ldots,
$$
\ifTrans
and we will have exactly the same calculations as for $ P_2 $, with only
the change of the cosines into sines; so we will have
\fi
$$
\begin{array}{p{0.3em}p{0.2em}p{18em}}
$p$&$=$&$a\sin(g+2)\alpha+c\sin(g-2)\alpha,$\\
$p_1$&$=$&$d\sin(g+4)\alpha+k\sin(g-4)\alpha,$\\
$p_2$&$=$&$l\sin(g+6)\alpha+n\sin(g+2)\alpha+q\sin(g-2)\alpha$\\
&&$+s\sin(g-6)\alpha,$\\
\end{array}
$$
$$
.\ .\ .\ .\ .\ .\ .\ .\ .\ .\ .\ .\ .\ .\ .\ .\ .\ .\ .\ .\ .\ .\ .\ .\ .\ .\ .\ .\ .\ .\ .\ .\ .\ .
$$
\ifTrans
and $a$, $c$, $d$, $k$ have the same values as in the expression of
$P_2$; so we still have for the constant
\fi
$$
\begin{array}{p{0.4em}p{0.2em}p{16em}}
$R'$&$=$&$g^2+\tfrac{1}{2(g^2-1)}h^4
+\tfrac{5g^2+7}{32(g^2-1)^3(g^2-4)}h^8$\\
&&$+\tfrac{9g^6+22g^4-203g^2-116}{64(g^2-1)^5(g^2-4)^2(g^2-9)}h^{12}
+\ldots,$
\end{array}
$$
\ifTrans
and this development must be stopped at the same term as in $R$;
then, although the first terms of $R$ and $R'$ are the same, these two
constants are not equal, and the two series separate from the term from
which we are forced to replace $g$ by its particular value.

{\bf 14}. When we have obtained the value of $P_2$ for an odd value
of $g$, it is easy to deduce that of $P_1$ for the same value of $g$.
Indeed, if $P_2$ is given by the equation
\fi
\begin{equation}
\tag{$n$}\label{eq:4n}
\tfrac{d^2P}{d\alpha^2}+\left[R(g^2,h^2)-2h^2\cos2\alpha\right]P=0,
\end{equation}
\ifTrans
by changing $h^2$ to $-h^2$ and $\alpha$ to $\tfrac{\pi}{2}-\alpha$, we
will have a periodic function which will satisfy the equation
\fi
\begin{equation}
\tag{$p$}\label{eq:4p}
\tfrac{d^2P}{d\alpha^2}+\left[R(g^2,-h^2)-2h^2\cos2\alpha\right]P=0,
\end{equation}
\ifTrans
same form as \eqref{eq:4m}, and if $g$ is odd, the cosines of $P_2$
change into sines; so we have the expression of $P_1$, and moreover we
see that we obtain the constant $R'$, which suits $P_1$, by changing in 
$R$ $h^2$ to $-h^2$.

According to this, for $g=1$ we have
\fi
$$
\begin{array}{p{0.4em}p{0.2em}p{20.5em}}
$P_1$&$=$&$\sin\alpha-\tfrac{h^2}{8}\sin3\alpha
+h^4\left(\tfrac{1}{192}\sin5\alpha+\tfrac{1}{64}\sin3\alpha\right)$\\
&&$-h^6\left(\tfrac{1}{9216}\sin7\alpha
+\tfrac{1}{1152}\sin5\alpha+\tfrac{1}{1536}\sin3\alpha\right)$\\
&&$+h^8\left(\tfrac{1}{737280}\sin9\alpha
+\tfrac{1}{49152}\sin7\alpha\right.$\\
&&$\qquad\quad+\left.\tfrac{1}{24576}\sin5\alpha
-\tfrac{11}{36864}\sin3\alpha\right)+\ldots;$
\end{array}
$$
$$
R'=1-h^2-\tfrac{1}{8}h^4+\tfrac{1}{64}h^6-\tfrac{1}{1536}h^8-\tfrac{11}{36864}h^{10}+\ldots.
$$
\ifTrans
And for $g=3$ we have
\fi
$$
\begin{array}{p{0.4em}p{0.2em}p{20.5em}}
$P_1$&$=$&$\sin3\alpha
+h^2\left(-\tfrac{1}{16}\sin5\alpha+\tfrac{1}{12}\sin\alpha\right)$\\
&&$+h^4\left(\tfrac{1}{640}\sin7\alpha+\tfrac{1}{64}\sin\alpha\right)$\\
&&$+h^6\left(\tfrac{-1}{46080}\sin9\alpha
-\tfrac{7}{20480}\sin5\alpha
+\tfrac{1}{768}\sin\alpha\right)$\\
&&$+h^8\left(\tfrac{1}{2^{14}\cdot 3^2\cdot 5\cdot7}\sin11\alpha
-\tfrac{17}{2^{5}\cdot3^2\cdot5}\sin7\alpha\right.$\\
&&$\qquad\quad-\left.\tfrac{1}{2^{14}}\sin5\alpha
+\tfrac{1}{2^{13}}\cos\alpha\right)+\ldots;$
\end{array}
$$

$$
R'=9+\tfrac{1}{16}h^4-\tfrac{1}{64}h^6+\tfrac{59}{61440}h^8
+\tfrac{3}{16384}h^{10}+\ldots.
$$
\ifTrans
If g is even and $P_2$ is given by the equation \eqref{eq:4n}, changing 
$h^2$ to $-h^2$ and $\alpha$ to $\tfrac{\pi}{2}-\alpha$ in $P_2$, we will
have a function $P$ which will satisfy the equation \eqref{eq:4p}; but
the cosines remain cosines in this change; so the new expression still
belongs to $P_2$, and we conclude
\fi
$$
R(g^2,-h^2)=R(g^2,h^2).
$$
\ifTrans
The same reasoning is applicable to $P_1$; therefore, if $g$ is even, 
$P_1$ does not change when we replace $\alpha$ by $\tfrac{\pi}{2}-\alpha$
and $h^2$ by $-h^2$, and $R'$ contains only fourth powers of $h$.

By a special calculation, we find, for $g=2$;
\fi
$$
\begin{array}{p{0.4em}p{0.2em}p{18em}}
$P_1$&$=$&$\sin2\alpha
-\tfrac{h^2}{12}\sin4\alpha
+\tfrac{h^4}{384}\sin6\alpha$\\
&&$+h^6\left(\tfrac{-1}{23040}\sin8\alpha
+\tfrac{5}{13824}\sin4\alpha\right)$\\
&&$+h^8\left(\tfrac{1}{2211840}\sin10\alpha
+\tfrac{37}{2209140}\sin6\alpha\right)$\\
&&$+h^{10}\left(\tfrac{-1}{309657600}\sin12\alpha
+\tfrac{11}{33177600}\sin8\alpha\right.$\\
&&$\qquad\quad-\left.\tfrac{289}{79626240}\sin4\alpha\right)+\ldots;$
\end{array}
$$
$$
R'=4-\tfrac{1}{12}h^4+\tfrac{5}{13824}h^8
-\tfrac{289}{79626240}h^{12}+\ldots
$$
\ifTrans
For $g=4$, we find
\fi
$$
\begin{array}{p{0.4em}p{0.2em}p{22em}}
$P_1$&$=$&$\sin4\alpha
+h^2\left(-\tfrac{1}{20}\sin6\alpha+\tfrac{1}{12}\sin2\alpha\right)
+\tfrac{h^4}{960}\sin8\alpha$\\
&&$-h^6\left(\tfrac{1}{80640}\sin10\alpha
+\tfrac{13}{96000}\sin6\alpha
-\tfrac{1}{4320}\sin2\alpha\right)$\\
&&$+h^8\left(\tfrac{1}{10321920}\sin12\alpha
+\tfrac{23}{6048000}\sin8\alpha\right)$\\
&&$+h^{10}\left(\tfrac{-1}{1857945600}\sin14\alpha
-\tfrac{53}{1032192000}\sin10\alpha\right.$\\
&&$\qquad\quad+\left.\tfrac{293}{2419200000}\sin6\alpha
+\tfrac{397}{124416000}\sin2\alpha\right)+\ldots;$
\end{array}
$$
$$
R'=16+\tfrac{1}{30}h^4-\tfrac{317}{864000}h^8
+\tfrac{4507}{1360800000}h^{12}+\ldots
$$

\ifTrans
\begin{center}
\emph{On the functions $Q$ which must be associated with $P_1$ and $P_2$.}
\end{center}

{\bf 15}. We go from the equation
\fi
\begin{equation}
\tag{1}\label{eq:51}
\tfrac{d^2P}{d\alpha^2}+\left(R-2h^2\cos2\alpha\right)P=0
\end{equation}
\ifTrans
to the one giving $Q$
\fi
\begin{equation}
\tag{2}\label{eq:52}
\tfrac{d^2Q}{d\beta^2}-\left[R-2h^2E(2\beta)\right]Q=0,
\end{equation}
\ifTrans
by changing $\alpha$ to $\beta i$ and $P$ to $Q$, and $R$ has the same
value in both; so if in the values of $P_1$ and $P_2$ we change $\alpha$
to $\beta i$, we will have solutions of the equation \eqref{eq:52}. Now
it is easy to understand that if the membrane folds in a vibratory
movement along the lines $FA$ and $F'A'$, included between the foci and
the vertices of the major axis, and which have the equations $\alpha=0$
and $\alpha=\pi$, it must also bend all along the line $FF'$ led between
the foci and which is given by $\beta=0$. So
\fi
$$
P_1=\sin g\alpha+h^2[a\sin(g+2)\alpha+b\sin(g-2)\alpha]+\ldots,
$$
\ifTrans
which becomes zero for $\alpha=0$ and $\alpha=\pi$, must effectively be
associated with the function
\fi
$$
Q_1=A\left\{\tfrac{e^{\beta g}-e^{-\beta g}}{2}
+h^2\left[a\tfrac{e^{\beta(g+2)}-e^{-\beta(g+2)}}{2}
+b\tfrac{e^{\beta(g+2)}-e^{-\beta(g+2)}}{2}
\right]+\ldots\right\},
$$
\ifTrans
which is deduced from $P_1$ by changing $\alpha$ to $\beta i$, and which
becomes zero for $\beta=0$. Likewise
\fi
$$
P_2=\cos g\alpha+h^2[a\cos(g+2)\alpha+b\cos(g-2)\alpha]+\ldots,
$$
\ifTrans
which is maximum or minimum for $\alpha=0$, must associate with
\fi
$$
Q_2=A\left\{E(\beta g)
+h^2\left[aE(\overline{g+2}\,\beta)
+bE(\overline{g-2}\,\beta)\right]\right\}+\ldots,
$$
\ifTrans
which is maximum or minimum for $\beta=0$.

However the expressions of $P_1$ and $P_2$ could be convergent, without
those of $Q_1$ and $Q_2$ being so, for all the values that $\beta$ can
take inside the membrane; but, for the moment, we want to point out the
characters of the functions $Q_1$ and $Q_2$ rather than to give a means
of calculating them.
\fi
\ifTrans
\begin{center}
\emph{On the nodal lines.}
\end{center}

{\bf 16}. We can already make some reflections on the nature of the nodal lines of an elliptical membrane. We have seen that in a
simple vibratory movement the displacement of a point of the membrane is given by the formula
\fi
$$
w= PQ \sin2\lambda mt,
$$
\ifTrans
where $P$ and $Q$ satisfy the two equations \eqref{eq:51} and 
\eqref{eq:52} of the previous number, and where $\lambda$ is determined
by the fixity of the contour. And $P$ being a function of $\alpha$ with
period $ 2\pi$, we can only take for it $P_1$ and $P_2$, so that $Q_1$
and $Q_2$ being the corresponding values of $Q$, we have two kinds of
solutions given by the formulas
\fi
$$
w=P_1Q_1\sin2\lambda_1\,mt,
$$
$$
w=P_2Q_2\sin2\lambda_2\,mt,
$$
\ifTrans
and the nodal lines have for equations, in the first kind,
\fi
$$
P_1=0,\quad Q_1=0,
$$
\ifTrans
and, in the second kind,
\fi
$$
P_2=0,\quad Q_2=0.
$$
\ifTrans
In the first kind, the major axis is a nodal line; in the second kind, there is maximum or minimum vibration.

The equations $Q_1=0$ and $Q_2=0$ give ellipses which have the same foci
as the contour of the membrane. The equations $P_1=0$ and $P_2=0$
determine the asymptotes of the hyperbolic nodal lines which still have
the same foci; and the integer $g$ which enters $P_1$ and $P_2$ indicates
how many times these functions become zero from $0$ to $\pi$, i.e.
the number of hyperbolic nodal lines, in designating by \emph{hyperbolic
nodal line} the two branches of a hyperbola terminated at the major axis
which have the same asymptote. In this view, a hyperbola is counted for
two of these lines; but if the
\fi
\ifTrans
major axis or the minor axis are motionless, they are only counted for a
single hyperbolic line.

If g is zero, the movement can only be of the second kind, and there is
no hyperbolic line.

If $g=1$, there is a hyperbolic nodal line only the major axis in the
first kind and only the minor axis in the second kind.

If $g=2$, in the first kind we have for these lines the small and the
long axis, and in the second kind a hyperbola.

If $g=3$, we have for these nodal lines a hyperbola and either the major
axis or the minor axis, depending on whether the movement is of the first
or of the second kind. And so on.

When the series found for $P_1$ and $P_2$ will be rapidly convergent, as
we know exactly the number of roots understood from $0$ to 
$\tfrac{\pi}{2}$, our formulas will be very convenient, and it will be
easy to separate these roots by substitutions and to calculate them with
the approximation desired by experience.

But if $h$, which is proportional to the eccentricity of the membrane
and the pitch of the sound, is large enough, these series are no longer
convergent or are too small to be of convenient use; we can no longer
even use the expansions of $R$ and $R'$ according to the powers of $h$
and we are forced to resort to other methods.

\begin{center}
\emph{Developments of the functions $P_1$ and $P_2$ according to powers
of $\sin\alpha$\\ and $\cos\alpha$}.
\end{center}

{\bf 17}. If we let
\fi
$$
\nu=\cos\alpha
$$
\ifTrans
and we take $\nu$ to be variable, the equation
\fi
\begin{equation}
\tag{1}\label{eq:71}
\tfrac{d^2P}{d\alpha^2}+\left[R(g^2,h^2)-2h^2\cos2\alpha\right]P=0
\end{equation}
\ifTrans
turns into the following :
\fi
\begin{equation}
\tag{$a$}\label{eq:7a}
\tfrac{d^2P}{d\nu^2}(1-\nu^2)-\tfrac{dP}{d\nu}\nu
+\left[R(g^2,h^2)+2h^2-4h^2\nu^2\right]P=0
\end{equation}
\ifTrans
and if we let
\fi
$$\nu'=\sin\alpha$$
\ifTrans
be variable, it changes to this other :
\fi
\begin{equation}
\tag{$b$}\label{eq:7b}
\tfrac{d^2P}{d\nu'^2}(1-\nu'^2)-\tfrac{dP}{d\nu'}\nu'
+\left[R(g^2,h^2)-2h^2+4h^2\nu'^2\right]P=0.
\end{equation}

\ifTrans
Suppose first $g$ is even: the constant $R(g^2,h^2)$, as we have seen,
depends only on the even powers of $h^2$; so we go from the equation 
\eqref{eq:7a} to the equation \eqref{eq:7b} by changing $\nu$ to $\nu'$
and $h^2$ to $-h^2$, and we conclude this remarkable property that, when
$g$ is even, $P_1$ and $P_2$ are functions of $\nu$ which remain
constant when we change $\nu$ to $\nu'$ and $h^2$ in $-h^2$.\footnote{Trans.\ note: The $h$ appearing in the original has been corrected to the $h^2$ appearing here.}

Suppose $g$ is odd: if $P_1$ is a solution of equation \eqref{eq:71}, we
know that the value of $P_2$ corresponding to the same value of $g$ is a
solution of the same equation in which we replace $R(g^2,h^2)$ by 
$R(g^2,-h^2)$; then $P_1$ also satisfies the two equations \eqref{eq:7a}
and \eqref{eq:7b}, and $P_2$ to the following two :
\fi
\begin{equation}
\tag{$c$}\label{eq:7c}
\tfrac{d^2P}{d\nu^2}(1-\nu^2)-\tfrac{dP}{d\nu}\nu
+\left[R(g^2,-h^2)+2h^2-4h^2\nu^2\right]P=0.
\end{equation}
\begin{equation}
\tag{$d$}\label{eq:7d}
\tfrac{d^2P}{d\nu'^2}(1-\nu'^2)-\tfrac{dP}{d\nu'}\nu'
+\left[R(g^2,-h^2)-2h^2+4h^2\nu'^2\right]P=0.
\end{equation}
\ifTrans
We go from \eqref{eq:7a} to \eqref{eq:7d} or from \eqref{eq:7c} to 
\eqref{eq:7b} by changing $\nu$ to $\nu'$ and $h^2$ in $-h^2$; therefore,
by the same changes, we pass from the expression of $P_1$ to that of 
$P_2$, or conversely from that of $P_2$ to that of $P_1$.

Let us consider equation \eqref{eq:7a} and let, to simplify the writing,
\fi
$$
R(g^2,h^2)+2h^2=m.
$$
\ifTrans
The general solution of this equation is the sum of two particular
solutions, one $\Pi_2$ even in $\nu$ and the other $\Pi_1$ odd. For the
function $\Pi_2$, we let
\fi
\begin{equation}
\tag{$e$}\label{eq:7e}
\Pi_2=k_0+k_1\nu^2+k_2\nu^4+k_3\nu^6+\ldots
+k_{s-1}\nu^{2s-2}+k_{s}\nu^{2s}+k_{s+1}\nu^{2s+2}+\ldots
\end{equation}
\ifTrans
and we will determine, by substituting in \eqref{eq:7a}, the coefficients
\fi
$$
k_1=-\tfrac{mk_0}{2},\quad k_2=\tfrac{m(m-4)+8h^2}{2\cdot3\cdot4}k_0,
$$
$$
k_3=\tfrac{-m(m-4)(m-16)-56h^2m+128h^2}{2\cdot3\cdot4\cdot5\cdot6}k_0,\ldots.
$$

\ifTrans
By equating to zero the coefficient of $\nu^{2s}$ in the result of the
substitution, we obtain the formula
\fi
\begin{equation}
\tag{2}\label{eq:72}
k_{s+1}=\tfrac{(4s^2-m)k_s+4h^2k_{s-1}}{(2s+1)(2s+2)},
\end{equation}
\ifTrans
which shows how each term is deduced from the previous two.

For the function $\Pi_1$, let
\fi
\begin{equation}
\tag{$f$}\label{eq:7f}
\left\{
\begin{array}{p{0.5em}p{0.2em}p{18em}}
$\Pi_1$&$=$&$a_1\nu+a_2\nu^3+a_3\nu^5+a_4\nu^7+\ldots$\\
&&$+a_{s-1}\nu^{2s-3}+a_s\nu^{2s-1}+a_{s+1}\nu^{2s+1}+\ldots,$
\end{array}
\right.
\end{equation}
\ifTrans
and we will have for the coefficients of the first terms\footnote{Trans.\ note: corrected the original by including a
multiplicative factor of $a_1$ missing in the second displayed equation
below.}
\fi
$$
a_2=-\tfrac{m-1}{2\cdot3}a_1,\quad a_3=\tfrac{(m-1)(m-9)+24h^2}{2\cdot3\cdot4\cdot5}a_1,
$$
$$
a_4=\tfrac{-(m-1)(m-9)(m-25)+(170-26m)4h^2}{2\cdot3\cdot4\cdot5\cdot6\cdot7}a_1,\ldots,
$$
\ifTrans
and each term is deduced from the two preceding ones by 
the relation
\fi
\begin{equation}
\tag{3}\label{eq:73}
a_{s+1}=\tfrac{[(2s-1)^2-m]a_s+4h^2s_{s-1}}{2s(2s+1)},\ldots.
\end{equation}

\ifTrans
Suppose the constant $R$ chosen so that $P_1$ satisfies the equation 
\eqref{eq:7a}; $P_1$ is zero or maximum for $\alpha=\tfrac{\pi}{2}$, and
behaves in this property like $\sin g\alpha$, to which it reduces, to
within a constant factor, for $h=0$; therefore it is zero if $g$ is even,
and it is maximum if $g$ is odd. However, for $\nu=0$ where 
$\alpha=\tfrac{\pi}{2}$, $\Pi_1$ is zero and $\Pi_2$ is maximum; so if 
$g$ is even $P_1$ is equal to $\Pi_1$, and if $g$ is odd $P_1$ is equal
to $\Pi_2$, except for a constant factor. Imagine, on the contrary, $R$
is chosen
\fi
\ifTrans
so that $P_2$ satisfies the equation \eqref{eq:7a}, and we also see that
$P_2$ is zero or maximum for $\alpha=\tfrac{\pi}{2}$, depending on
whether $g$ is odd or even; and we conclude that $P_2$ is equal, except
for a factor, to $\Pi_1$ if $g$ is odd, and to $\Pi_2$ if $g$ is even.

{\bf 18}. Let us come to the equation \eqref{eq:7b}, and let
\fi
$$
R-2h^2=m';
$$
\ifTrans
the general solution of this equation is also the sum of two particular
solutions, one of which is even and the other odd in $\nu'$. If $R$ is
chosen so that $P_2$ is solution of this equation, then $P_2$, which is
maximum for $\nu'=0$, merges with the particular solution which enjoys
this property, and we have
\fi
\begin{equation}
\tag{$e'$}\label{eq:7e2}
P_2=k'_0+k'_1\nu'^2+k'_2\nu'^4+\ldots+k'_{s-1}\nu'^{2s-2}+k'_s\nu'^{2s}+k'_{s+1}\nu'^{2s+2}+\ldots,
\end{equation}
\ifTrans
an expression where $k'_0$, $k'_1$, $k'_2$,\ldots are deduced from $k_0$,
$k_1$, $k_2$,\ldots, by changing $h^2$ to $-h^2$. So we have\footnote{Trans.\ note: I have corrected the subscript of the second equation to read $k'_2$ rather than $k'_1$.}
\fi
$$
k'_1=-\tfrac{m'k'_0}{1\cdot2},\quad
k'_2=\tfrac{m'(m'-4)-8h^2}{1\cdot2\cdot3\cdot4}k'_0,\ldots,
$$
\ifTrans
and the rule which links the coefficients of three consecutive terms
together is
\fi
\begin{equation}
\tag{4}\label{eq:74}
k'_{s+1}\tfrac{(4s^2-m')k'_s-4h^2k'_{s-1}}{(2s+1)(2s+2)}.
\end{equation}

\ifTrans
If, on the contrary, $R$ is chosen so that $P_1$ is a solution of
equation \eqref{eq:71}, $P_1$, having to become zero for $\nu'=0$, merges
with the odd solution in $\nu'$ of \eqref{eq:7b}, and we have
\fi
\begin{equation}
\tag{$f'$}\label{eq:7f2}
P_1=a'_1\nu'+a'_2\nu'^3+a'_3\nu'^5+\ldots+a'_{s-1}\nu'^{2s-3}+a'_s\nu'^{2s-1}+a'_{s+1}\nu'^{2s+1}+\ldots
\end{equation}
\ifTrans
$a'_1$, $a'_2$, $a'_3$,\ldots being quantities deduced from $a_1$, $a_2$,
$a_3$,\ldots by changing $h^2$ to $-h^2$; so we have first
\fi
$$
a'_2=-\tfrac{m'-1}{2\cdot3}a'_1,\quad
a'_3=\tfrac{(m'-1)(m'-9)-24h^2}{2\cdot3\cdot4\cdot5}a'_1,\ldots,
$$
\ifTrans
and then the general formula
\fi
\begin{equation}
\tag{5}\label{eq:75}
a'_{s+1}=\tfrac{[(2s-1)^2-m']a'_s-4h^2a'_{s-1}}{2s(2s+1)}.
\end{equation}

\ifTrans
Here it is very important to note that if we replace $R$ by an
arbitrary number, the function \eqref{eq:7e2} or the function
\eqref{eq:7f2} is not equal to any of the two functions \eqref{eq:7e} and
\eqref{eq:7f}; because the function \eqref{eq:7f2}, for example, which is
zero for $\nu'=0$ or $\alpha=0$, is neither zero nor maximum for 
$\alpha=\tfrac{\pi}{2}$ and cannot therefore be confused with either of
the two functions \eqref{eq:7e} and \eqref{eq:7f}. It is only if $R$ is
determined so that equation \eqref{eq:71} has a solution of period
$2\pi$, that, this solution having to be zero or maximum for $\alpha=0$
and $\alpha=\tfrac{\pi}{2}$, one of the two expressions \eqref{eq:7e2}
and \eqref{eq:7f2}, which is equal to it, is identical, except for a
factor, to one of the two expressions \eqref{eq:7e} and \eqref{eq:7f}.

For example, suppose $g$ even, and, consequently, $P_2$ is equal to the
product of $\Pi_2$ by a constant $A$, if $R$ has been taken correctly. By
considering the values of these functions for $\alpha=45^{\circ}$, 
$30^{\circ}$, $60^{\circ}$, we obtain the formulas
\fi
$$
A\left(k_0+\tfrac{k_1}{2}+\tfrac{k_2}{4}+\tfrac{k_3}{16}+\cdots\right)=
k'_0+\tfrac{k'_1}{2}+\tfrac{k'_2}{4}+\ldots,
$$
$$
A\left(k_0+k_1\tfrac{3}{4}+k_2\tfrac{9}{16}+\cdots\right)=
k'_0+\tfrac{k'_1}{4}+\tfrac{k'_2}{16}+\ldots,
$$
$$
\tfrac{1}{A}\left(k'_0+k'_1\tfrac{3}{4}+k'_2\tfrac{9}{16}+\cdots\right)=
k_0+\tfrac{k_1}{4}+\tfrac{k_2}{16}+\ldots,
$$
\ifTrans
each of which can determine the factor $A$.

We have so far accepted that $R$ was known; but if it is not, by
eliminating $A$ between two of these formulas, we will obtain an equation
whose two members will be the products of two series, which will contain only the unknown $R$ and may be used to determine it.

It is extremely easy to recognize that the series \eqref{eq:7e},
\eqref{eq:7f}, \eqref{eq:7e2}, \eqref{eq:7f2} are convergent, as long as
$\nu$ or $\nu'$ is $<1 $. Or more generally the series
\fi
$$
k_0+k_1x+k_2x^2+\ldots+k_nx^n+k_{n+1}x^{n+1}+\ldots,
$$
\ifTrans
in which $x$ is $<1 $, and of which three consecutive coefficients are
linked by the relation
\fi
\begin{equation}
\tag{7}\label{eq:77}
k_{s+1}=A_sk_s+a_sk_{s-1};
\end{equation}
\ifTrans
moreover the limit of $A_s$, when $s$ grows indefinitely, is less than
unity or it is at most equal, and the limit of $a_s$ is zero; then the
series is convergent.

Indeed, the limit of the ratio $\tfrac{k_{s+1}}{k_s}$ is equal to the
limit $\tau$ of $A_s$ when $s$ grows indefinitely; therefore the limit of
the relation of a term to the previous in the series is $\tau x$, a number
$<1$, and it is convergent. Now the relations \eqref{eq:72},
\eqref{eq:73}, \eqref{eq:74}, \eqref{eq:75}, satisfying the same
conditions as the relation \eqref{eq:77}, the series are convergent as
long as $\nu$ and $\nu'$ are $<1$, i.e. whatever $\alpha$ may be.

On the contrary, if $\nu$ and $\nu'$ were $>1$, these series would be
divergent.

However to have well convergent series, we will prefer the formulas
\eqref{eq:7e2} and \eqref{eq:7f2} when $\alpha$ will be between $0$ and 
$45$ degrees, and the formulas \eqref{eq:7e} and \eqref{eq:7f} when
$\alpha$ will be between $45$ and $90$ degrees.

{\bf 19}. It is easy to see that, from a sufficiently distant term, all
the following terms have the same sign in the four series \eqref{eq:7e},
\eqref{eq:7f}, \eqref{eq:7e2}, \eqref{eq:7f2} but we will also study the
number of variations of the latter two. If we first assume $h=0$, $R$ and
$m'$ are reduced to $g^2$, and the two series \eqref{eq:7e2} and 
\eqref{eq:7f2} to
\fi
$$
\def\arraystretch{1.8}
\begin{array}{p{0.4em}p{0.2em}p{22em}}
$P_2$&$=$&$k'_0\left[
1-\tfrac{g^2}{2}\nu'^2
+\tfrac{g^2(g^2-4)}{2\cdot3\cdot4}\nu'^4
-\tfrac{g^2(g^2-4)(g^2-16)}{2\cdot3\cdot4\cdot5\cdot6}\nu'^6+\ldots
\right],$\\
$P_1$&$=$&$\tfrac{a'_1}{g}\left[
g\nu'-\tfrac{g(g-1)}{2\cdot3}\nu'^3
+\tfrac{g(g^2-1)(g^2-9)}{2\cdot3\cdot4\cdot5}\nu'^5-\ldots
\right],$
\end{array}
$$
\ifTrans
except for factors close to the values of $\cos g\alpha$ and 
$\sin g\alpha$.

If we make successively
\fi
$$
g=1, 2, 3,\ldots,
$$
\ifTrans
the factor in square brackets of $P_2$ becomes\footnote{Trans.\ note: I added a factor of $\nu'^4$ in the fourth equation, which seems to fit, but this needs to be verified correct. RMC: This looks correct to me now.}
\fi
$$
\def\arraystretch{1.4}
\begin{array}{p{2em}p{0.2em}p{26em}}
$\ \ \cos\alpha$&$=$&$1-\tfrac{1}{1\cdot2}\nu'^2
-\tfrac{3}{1\cdot2\cdot3\cdot4}\nu'^4
-\tfrac{3\cdot15}{1\cdot2\cdot3\cdot4\cdot5\cdot6}\nu'^6-\ldots,$\\
$\cos2\alpha$&$=$&$1-2\nu'^2$\\
$\cos3\alpha$&$=$&$1-\tfrac{3^2}{1\cdot3}\nu'^2
+\tfrac{3^2(3^2-2^2)}{1\cdot2\cdot3\cdot4}\nu'^4
+\tfrac{3^2(3^2-2^2)(4^2-3^2)}{1\cdot2\cdot3\cdot4\cdot5\cdot6}\nu'^6+\ldots,$\\
$\cos4\alpha$&$=$&$1-\tfrac{4^2}{1\cdot2}\nu'^2
+\tfrac{4^2(4^2-2^2)}{1\cdot2\cdot3\cdot4}\nu'^4,$\\
\end{array}
$$
$$
.\ .\ .\ .\ .\ .\ .\ .\ .\ .\ .\ .\ .\ .\ .\ .\ .\ .\ .\ .\ .\ .\ .\ .\ \qquad\qquad\qquad\qquad\qquad\qquad\qquad\quad
$$

\ifTrans
The series which gives $\cos\alpha$ has a single variation
and a single positive root in $\nu'$, $\nu'=\sin\tfrac{\pi}{2}$; 
$\cos2\alpha$ also only has a variation and a positive root 
$\nu'=\sin\tfrac{\pi}{4}$; $\cos3\alpha$ has two variations and two
positive roots, $\nu'=\sin\tfrac{\pi}{6}$ and $\sin\tfrac{\pi}{2}$; 
$\cos4\alpha$ has two variations and two positive roots, 
$\nu'=\sin\tfrac{\pi}{8}$ and $\sin\tfrac{3\pi}{8}$. And in general the
series which gives $\cos g\alpha$ by means of $\nu'$ has as many
variations as the equation $\cos g\alpha=0$ has positive roots in 
$\nu'$.

We also recognize that the series which expresses $\sin g\alpha$ has a
number of variations equal to the number of its roots in $\nu'$.

So when $h$ is zero, $P_1$ and $P_2$ expressed by $\nu'$ have the same
number of variations as positive roots, and we will show that this
property subsists for any value of $h$.

Suppose, for example, that it is $P_2$; we have between the coefficients
of three consecutive terms of
\fi
$$
P_2=k'_0+k'_1\nu'^2+k'_2\nu'^4+\ldots
$$
\ifTrans
the relation
\fi
$$
k'_{s+1}=\tfrac{(4s^2-m')k'_s-4h^2k'_{s-1}}{(2s+1)(2s+2)},
$$
\ifTrans
and imagine that we increase the quantity $h$. It follows from this
formula that for any value of $h$ (except zero), two consecutive
coefficients cannot become zero. Indeed if $k'_s$, and 
$k'_{s+1}$, were zero, $k'_{s-1}$
\fi
\ifTrans
so would be, then for the same reason $k'_{s-2}$, and so on, so that all
of the terms in the series would become zero.

Secondly, if the coefficient of one of the terms becomes zero for a
certain value of $h$, the coefficients of the two terms which surround it
are of opposite sign, as we see by the same formula.

It follows from this that, while $ h $ increases, $ P_2 $ cannot acquire or
lose any variation, and that it consequently has the same number as for 
$h=0$. But, as we have shown (n$^{\circ}${\bf 10}), $P_2$ Always becomes zero
the same number of times $\alpha=0$ to $\alpha=\tfrac{\pi}{2}$,
regardless of $h$; so finally the equation
\fi
$$
P_2(\nu')=0
$$
\ifTrans
has precisely as many real, positive and $<1$ roots as it has variations.

{\bf 20}. This property separates the roots of this equation. First consider an algebraic equation
\fi
$$
f(x)=0,
$$
\ifTrans
which has as many positive roots as variations; let us 
form the sequence of
derivatives of $f(x)$
\fi
\begin{equation}
\tag{A}\label{eq:7A}
f(x), f'(x), f''(x),\ldots,
\end{equation}
\ifTrans
which for $x=0$ has the same signs as the series of coefficients of 
$f(x)$; it is easy to prove that, while one increases $x$, it is
impossible that this series never gains variations; but that when $f(x)$
becomes zero, a variation from the first to the second term is lost. So
if the equation $f(x)=0$ has as many variations as there are positive
roots, like the sequence \eqref{eq:7A} for $x=0$ has this number of
variations, that when $x$ increases, it loses one each time $ f (x) $ becomes
zero, and it cannot gain one, it can only lose one when $x$ passes
through a root of the equation $f(x)=0$, and counting the number of
variations of the sequence \eqref{eq:7A} for $x=a$ and $x=b$, and making
the difference, we have precisely the number of roots between $a$ and $b$.
\fi

\ifTrans
All these reasonings are applicable to the equation $P_2(\nu')=0$, formed
of a sequence of an infinite number of terms, and we will have to examine
the infinite sequence
\fi
\begin{equation}
\tag{B}\label{eq:7B}
P_2(\nu'),\quad \tfrac{dP_2}{d\nu'},\quad \tfrac{d^2P_2}{d\nu'^2},\quad \tfrac{d^3P_2}{d\nu'^3},\ldots,
\end{equation}
\ifTrans
the first two will be obtained by the sequence
\fi
$$
P_2=k'_0+k'_1\nu'^2+k'_2\nu'^4+\ldots,
$$
$$
\tfrac{dP_2}{d\nu'}=k'_12\nu'+k'_24\nu'^3+\ldots;
$$
\ifTrans
then we will calculate the following derivatives using the equation
\fi
$$
\tfrac{d^2P}{d\nu'^2}(1-\nu'^2)-\tfrac{dP}{d\nu'}\nu'+
(R-2h^2+4h^2\nu'^2)P=0
$$
\ifTrans
and those that we deduce by differentiation.

The infinite number of terms in the sequence \eqref{eq:7B} does not offer
any trouble; because we know how many roots $P_2(\nu')=0$ has between 
$0$ and $1$; it has as many as $\cos g\alpha=0$ between $0$ and 
$\tfrac{\pi}{2}$: it is $\tfrac{g}{2}$ or $\tfrac{g + 1}{2}$, depending
on whether $g$ is even or odd; $P_2(\nu')$ will have as many variations;
we will therefore calculate only a number $n$ of terms of the series
which gives $P_2$, sufficient to count all the variations there, and it
follows from the first principles of algebra that the sequence 
\eqref{eq:7B} will have no variations beyond its first $n$ terms when we
give $\nu'$ a positive value.

{\bf 21}. Everything we just said about $P_2$ can be repeated for $P_1$.
We can now get a more precise idea of the constant $R$. This quantity,
for $h=0$, is reduced to the square of an integer: it is therefore
positive when $h$ is very small; but we will demonstrate that the
constant $R$ relative to $P_2$ is not only positive, but also greater
than $2h^2$, whatever $h$.

We have just examined the changes in the number of variations of the
sequence \eqref{eq:7B} when we vary $\nu'$ from $0$ to $1$; but as
\fi
\ifTrans
$P_1$ and $P_2$ are one an odd function in $\nu'$ and the other an even
function, it follows that the series \eqref{eq:7B} has the same property
between $-1$ and $0$ it has between $0$ and $1$, and consequently if we
increase $\nu'$ from $-1$ to $+1$, a variation is lost only each time 
$\nu'$ goes through a root of $P_2(\nu')=0$. This happens for $P_1(\nu')$.

$P_2$ is given by the equation \footnote{Trans.\ note: corrected the missing '$ 2 $' in the second derivative.}
\fi
\begin{equation}
\tag{$n$}\label{eq:7n}
\tfrac{d^2P_2}{d\nu'^2}(1-\nu'^2)=\tfrac{dP_2}{d\nu'}\nu'+
(2h^2-4h^2\nu'^2-R)P_2;
\end{equation}
\ifTrans
let us make $\nu'=0$, we then have $\tfrac{dP_2}{d\nu'}=0$, and therefore
$P_2$ and $\tfrac{d^2P_2}{d\nu'^2}$ are of opposite sign; because if they
were of the same sign, letting $\nu'$ grow from a very small negative quantity
to a very small positive quantity, $\tfrac{dP_2}{d\nu'}$, becoming zero,
will pass from a sign opposite to that of $P_2$ to a similar sign, and
the series would lose two variations, while it should not lose any. So
the coefficient of $P_2$ is negative, and we have
\fi
$$
2h^2-R<0,
$$
\ifTrans
or $R>2h^2$.

Let us denote, as we have already done, by $R'$ the constant $R$ when it
belongs to the function $P_1$, and we will demonstrate that it is
$>-2h^2$ whenever the integer $g$ is $>1$. Indeed, $\tfrac{dP_1}{d\nu'}$
necessarily becomes zero for a value of $\nu'$ ranging from $0$ to $1$,
and, for this value, $P_1$, and $\tfrac{d^2P_1}{d\nu'^2}$ are of opposite
sign, and since $P_1$ satisfies the equation \eqref{eq:7n}, when we
replace $R$ by $R'$, we have
\fi
$$
2h^2-4h^2\nu'^2-R'<0
$$
\ifTrans
or
\fi
$$
R'>2h^2(1-2\nu'^2),
$$
\ifTrans
and even more so $R'$ is $>2h^2$.

It is good to notice that the equation \eqref{eq:7n} ceases to be
applicable to the limit $\nu'=1$, for this reason that $\nu'$, being the
sine of a determined angle, cannot take values greater than unity. And
\fi
\ifTrans
indeed, for $\nu'=1$, $P$ or $\tfrac{dP}{d\alpha}$ is zero. Suppose it is
$P_1$, it follows from this equation that $\tfrac{dP}{d\nu'}$ would be
zero, and therefore also $\tfrac{dP}{d\alpha}$, according to equality
\fi
$$
\tfrac{dP}{d\alpha}=\tfrac{dP}{d\nu'}\cos\alpha;
$$
\ifTrans
which is impossible.

\begin{center}
\emph{Differential equations which determine the function $Q$.}
\end{center}

{\bf 22}. We know that $Q$ is given by the equation
\fi
$$
\tfrac{d^2Q}{d\beta^2}-[R-2h^2E(2\beta)]=0
$$
\ifTrans
and if we let
\fi
$$
\rho=c\tfrac{e^{\beta}+e^{-\beta}}{2},\quad
\rho'=c\tfrac{e^{\beta}-e^{-\beta}}{2},
$$
\ifTrans
and we take $\rho$ and $\rho'$ for variables, we get the two equations
\fi
\begin{equation}
\tag{1}\label{eq:81}
\tfrac{d^2Q}{d\rho^2}(\rho^2-c^2)+\tfrac{dQ}{d\rho}\rho+
(4\lambda^2\rho^2-R-2h^2)Q=0,
\end{equation}
\begin{equation}
\tag{2}\label{eq:82}
\tfrac{d^2Q}{d\rho'^2}(\rho'^2-c^2)+\tfrac{dQ}{d\rho'}\rho'+
(4\lambda^2\rho'^2-R-2h^2)Q=0.
\end{equation}

\ifTrans
When we make $c=0$ in these equations, they are reduced to one, relating
to the circular membrane; the two semi-axes $\rho$ and $\rho'$ of any of
the ellipses confocal to the membrane changing into the radius $r$ of a
circle, we have
\fi
\begin{equation}
\tag{3}\label{eq:83}
\tfrac{d^2Q}{dr^2}r^2+\tfrac{dQ}{dr}r+
(4\lambda^2r^2-g^2)Q=0,
\end{equation}
\ifTrans
an equation found at n$^{\circ}${\bf 2}. We have seen that its general
solution is the sum of two particular solutions, one of which becomes
infinite for $r=0$. It must be explained how the solution of the circle
can be deduced from that of the ellipse.

For this, first suppose $\lambda$ is zero in the equations \eqref{eq:81},
\eqref{eq:82} and \eqref{eq:83};
\fi
\ifTrans
they become, by noting that $R$ is reduced to $g^2$ for the hypothesis 
$\lambda=0$, which gives $h=0$,
\fi
\begin{equation}
\tag{1$'$}\label{eq:812}
\tfrac{d^2Q}{d\rho^2}(\rho^2-c^2)+\tfrac{dQ}{d\rho}\rho-g^2Q=0,
\end{equation}
\begin{equation}
\tag{2$'$}\label{eq:822}
\tfrac{d^2Q}{d\rho'^2}(\rho'^2-c^2)+\tfrac{dQ}{d\rho'}\rho'-g^2Q=0,
\end{equation}
\begin{equation}
\tag{3$'$}\label{eq:832}
\tfrac{d^2Q}{dr^2}r^2+\tfrac{dQ}{dr}r-g^2Q=0.
\end{equation}

\ifTrans
The general integral of \eqref{eq:832} is
\fi
$$
Q=Ar^g+Br^{-g},
$$
\ifTrans
and becomes infinite for $r=0$, that is to say at the center of the
circle. But the integrals of the two equations \eqref{eq:812} and 
\eqref{eq:822} are
\fi
$$
Q=A\left(\rho+\sqrt{\rho^2-c^2}\right)^g+
\tfrac{B}{\left(\rho\,+\,\sqrt{\rho^2\,-\,c^2}\right)^g},
$$
$$
Q=A\left(\rho'+\sqrt{\rho'^2+c^2}\right)^g+
\tfrac{B}{\left(\rho'\,+\,\sqrt{\rho'^2\,+\,c^2}\right)^g},
$$
\ifTrans
and we see that the second part of their expression is not infinite for 
$\rho=c$ where $\rho'=0$, except in the case where $c$ is zero.

Thus the general solution of the equations \eqref{eq:81} and 
\eqref{eq:82}, when we make $\lambda=0$, contains two arbitrary constants
and does not become infinite for $\rho'=0$.

We have seen that the displacement of a point on the membrane is represented by
\fi
$$
w=Au\sin2\lambda mt,\quad u=PQ,
$$
\ifTrans
and for the same value of $g$, $P$ can become zero or be maximum for 
$\alpha=0$, and has two expressions, $P_1$ and $P_2$, which correspond to
two different values of the constant $R$, which only become identical for
$h=0$; from there, for $u$, two expressions,
\fi
$$
u=P_1Q_1,\quad u=P_2Q_2,
$$
\ifTrans
in which $Q_1$ is a value of $Q$ which becomes zero for $\beta=0$, and 
$Q_2$ a value which becomes maximum for this value of $\beta$. Let us see
what $Q_1$ and $Q_2$ reduce to when we make $\lambda=0$.

Let us say that
\fi
$$
Q=A\left(\rho'+\sqrt{\rho'^2+c^2}\right)^g+
B\left(\rho'+\sqrt{\rho'^2+c^2}\right)^{-g}
$$
\ifTrans
is zero for $\rho'=0$, and we will have $\tfrac{B}{A}=-c^{2g}$; so
\fi
$$
Q_1=A\left[\left(\rho'+\sqrt{\rho'^2+c^2}\right)^g-
\tfrac{c^{2g}}{\left(\rho'\,+\,\sqrt{\rho'^2\,+\,c^2}\right)^g}
\right];
$$
\ifTrans
and if we express that $\tfrac{dQ}{d\rho'}$ is zero for $\rho'=0$, we have
$\tfrac{B}{A}=c^{2g}$; so
\fi
$$
Q_2=A\left[\left(\rho'+\sqrt{\rho'^2+c^2}\right)^g+
\tfrac{c^{2g}}{\left(\rho'\,+\,\sqrt{\rho'^2\,+\,c^2}\right)^g}
\right].
$$

\ifTrans
Finally, if we make $c=0$, the two values of $Q_1$ and $Q_2$ for the same
value of the integer $g$ are identical.

These explanations were useful in helping to understand how the
theory of the circular membrane is enclosed in that of the elliptical
membrane; because it is obvious that what we have just found when 
$\lambda$ is zero, is true for any value of $\lambda$.

{\bf 25}. Let us return to equations \eqref{eq:81} and \eqref{eq:82}. If
we let
\fi
$$
\tfrac{\rho}{c}=u,
$$
\ifTrans
we will have, instead of the equation \eqref{eq:81},
\fi
$$
\tfrac{d^2Q}{du^2}(u^2-1)+\tfrac{dQ}{du}+(4h^2u^2-R-2h^2)Q=0,
$$
\ifTrans
and this equation is deduced from that which gives $P$ by means of $\nu$
by only the change from $P$ to $Q$ and from $\nu$ to $u$. Now we have
seen that $P_2$ is given by the series
\fi
$$
P_2=k_0+k_1\nu^2+k_2\nu^4+\ldots,
$$
\ifTrans
and $P_1$ by this other
\fi
$$
P_1=a_1\nu+a_2\nu^3+a_3\nu^5+\ldots,
$$
\ifTrans
$k_0$, $k_1$, $k_2$,\ldots, $a_1$, $a_2$, $a_3$,\ldots having the values
calculated at n$^{\circ}${\bf 17}; so the corresponding values of $Q$
are\footnote{Trans.\ note: I believe that the second equation should be for $Q_1$ rather than $Q$, a correction I have made.}
\fi
$$
Q_2=k_0+k_1\tfrac{\rho^2}{c^2}+k_2\tfrac{\rho^4}{c^4}+\ldots,
$$
$$
Q_1=a_1\tfrac{\rho}{c}+a_2\tfrac{\rho^3}{c^3}+a_3\tfrac{\rho^5}{c^5}+\ldots.
$$

\ifTrans
However, these series cannot be used; because they are never convergent
if $\rho$ is $>c$; which always takes place in our problem.

But let
\fi
$$
\tfrac{\rho'}{c'}=u',
$$
\ifTrans
instead of the equation \eqref{eq:82}
\fi
$$
\tfrac{d^2Q}{du'^2}(u'^2+1)+\tfrac{dQ}{du'}u'+
(4h^2u'^2-R+2h^2)Q=0,
$$
\ifTrans
which is deduced from the equation which gives $P$ by means of $\nu'$ 
(n$^{\circ}${\bf 17}), by changing $\nu'$ to $u'\sqrt{-1}$. So $Q_2$ and
$Q_1$ will be given by the formulas
\fi
$$
\def\arraystretch{1.6}
\begin{array}{p{0.5em}p{0.2em}p{14em}}
$Q_2$&$=$&$k'_0-k'_1\tfrac{\rho'^2}{c^2}+k'_2\tfrac{\rho'^4}{c^4}
-k'_3\tfrac{\rho'^6}{c^6}+\ldots,$\\
$Q_1$&$=$&$a'_1\tfrac{\rho'}{c}-a'_2\tfrac{\rho'^3}{c^3}+a'_3\tfrac{\rho'^5}{c^5}-\ldots.$
\end{array}
$$

\ifTrans
We know that the series which give $P_1$ and $P_2$ by means of $\nu'$ are
convergent as long as $\nu'$ is $<1$; the previous ones are deduced by
the change from $\nu'$ to $\tfrac{\rho'}{c}\sqrt{-1}$; it therefore
follows from the convergence circle theorem that these series are
convergent, as long as $\rho'$ is $<c$. These series will be very
convenient, if we consider a very-eccentric membrane so that $\rho'$ is
much smaller than $c$.
\fi

\ifTrans
When $\rho'$ is $>c$, we can usually get $Q$ with enough approximation as
follows. Let
\fi
$$
z=\tfrac{\rho\,+\,\sqrt{\rho^2\,-\,c^2}}{2}=\tfrac{ce^{\beta}}{2},
$$
\ifTrans
and taking $z$ to be variable, we will have
\fi
$$
z^2\tfrac{d^2Q}{dz^2}+z\frac{dQ}{dz}+\left[
4\lambda^2z^2\left(1+\tfrac{c^4}{16z^4}\right)-R
\right]Q=0.
$$
\ifTrans
However, if $\rho'=\sqrt{\rho^2-c^2}$ is $>c$, we will have
\fi
$$
\rho>c\sqrt{2},\quad\tfrac{c}{z}<2(\sqrt{2}-1),
$$
\ifTrans
$$
\tfrac{c^4}{16z^4}<17-12\sqrt{2}\quad\rm{or}\quad<0.03056\ldots.
$$
\fi

\ifTrans
So if in the previous differential equation, we reduce the factor of 
$4\lambda^2z^2$ to unity, the one that will result will make $Q$ known in
general with some approximation. Now the equation then takes the form
that we found for the circle
\fi
$$
z^2\tfrac{d^2Q}{dz^2}+z\tfrac{dQ}{dz}+(4\lambda^2z^2-R)Q=0,
$$
\ifTrans
and by setting $R=n^2$, we will have for approximate value the expression
\fi
$$
\begin{array}{p{0.4em}p{0.2em}p{20em}}
$Q$&$=$&$Cz^n\left[1-\tfrac{(\lambda z)^2}{1(n+1)}+
\tfrac{(\lambda z)^4}{1\cdot2(n+1)(n+2)}\right.$\\
&&$\qquad\quad\,-\left.\tfrac{(\lambda z)^6}{1\cdot2\cdot3(n+1)(n+2)(n+3)}
+\cdots\right],$
\end{array}
$$
\ifTrans
where $n$ however is no longer an integer as in the case of the circle,
but depends on $h$. Besides, we will give another way to develop $Q$
later.

\begin{center}
\emph{Developments of $Q_1$ and $Q_2$.}
\end{center}

{\bf 24}. The functions\footnote{Trans.\ note: the original reads ``The functions $Q$ and $Q_1$ satisfy\ldots'', which I have corrected.} $Q_2$ and $Q_1$ satisfy the equation
\fi
$$
\tfrac{d^2Q}{d\beta^2}=\left[R-h^2(e^{2\beta}+e^{-2\beta})\right]Q;
$$
\ifTrans
let
\fi
$$
R-h^2(e^{2\beta}+e^{-2\beta})=T,
$$
\ifTrans
and form the derivatives
\fi
$$
\tfrac{dT}{d\beta}=-2h^2(e^{2\beta}-e^{-2\beta}),\quad
\tfrac{d^2T}{d\beta^2}=-2^2h^2(e^{2\beta}+e^{-2\beta}),\ldots.
$$

\ifTrans
If we make $\beta=0$, all the derivatives of odd order become zero,
and designating the others by $A_2$, $A_4$,\ldots, let
\fi
$$
T=A_0,\quad\tfrac{d^2T}{d\beta^2}=-2^3h^2=A_2,\ldots,\quad
\tfrac{d^{2i}T}{d\beta^{2i}}=-2^{2i+1}h^2=A_{2i};
$$
\ifTrans
then develop $Q$ using the formula
\fi
$$
Q=Q_0+\left(\tfrac{dQ}{d\beta}\right)_0\beta+
\left(\tfrac{d^2Q}{d\beta^2}\right)_0\tfrac{\beta^2}{1\cdot2}+\ldots,
$$
\ifTrans
and form the derivatives using the formulas
\fi
$$
\def\arraystretch{1.6}
\begin{array}{p{1em}p{0.2em}p{11em}}
$\tfrac{d^2Q}{d\beta^2}$&$=$&$TQ,\quad
\tfrac{d^3Q}{d\beta^3}=Q\tfrac{dT}{d\beta}+T\tfrac{dQ}{d\beta},$\\
$\tfrac{d^4Q}{d\beta^4}$&$=$&$Q\tfrac{d^2T}{d\beta^2}+
2\tfrac{dQ}{d\beta}\tfrac{dT}{d\beta}+T\tfrac{d^2Q}{d\beta^2},$
\end{array}
$$
$$
.\ .\ .\ .\ .\ .\ .\ .\ .\ .\ .\ .\ .\ .\ .\ .\ .\ .\ .\ .\ .\ .\ .
$$

\ifTrans
Let's deal with $Q_1$ first; since it is zero for $\beta=0$, all the even
order derivatives of $Q$ become zero, and the odd order
derivatives have the values, by designating the first by $B$:
\fi
$$
\def\arraystretch{1.6}
\begin{array}{p{2.5em}p{0.2em}p{21em}}
$\ \left(\tfrac{dQ}{d\beta}\right)_0$&$=$&$B,\quad
\left(\tfrac{d^3Q}{d\beta^3}\right)_0=BA_0,\quad
\left(\tfrac{d^5Q}{d\beta^5}\right)_0=B(A_0^2+3A_2),$\\
$\left(\tfrac{d^7Q}{d\beta^7}\right)_0$&$=$&$B\left(A_0^3+
\tfrac{2\cdot3+4\cdot5}{1\cdot2}A_0A_2+
\tfrac{2\cdot3\cdot4\cdot5}{1\cdot2\cdot3\cdot4}A_4\right),$\\
$\left(\tfrac{d^9Q}{d\beta^9}\right)_0$&$=$&$B\left(A_0^4+
\tfrac{2\cdot3+4\cdot5+6\cdot7}{1\cdot2}A_0^2A_2+
\tfrac{2\cdot3\cdot4\cdot5+4\cdot5\cdot6\cdot7}{1\cdot2\cdot3\cdot4}A_0A_4\right.$\\
&&$\qquad\qquad\qquad\qquad\qquad\left.
\tfrac{2\cdot3\times 6\cdot7}{1\cdot2\times 1\cdot2}A_2^2+
\tfrac{2\cdot3\cdot4\cdot5\cdot6\cdot7}{1\cdot2\cdot3\cdot4\cdot5\cdot6}A_6\right).$
\end{array}
$$
$$
.\ .\ .\ .\ .\ .\ .\ .\ .\ .\ .\ .\ .\ .\ .\ .\ .\ .\ .\ .\ .\ .\ .\ .\ .\ .\ .\ .\ .\ .\ .\ .\ .\ .\ .\ .\ .\ .\ .\ .\ .\ .
$$
\ifTrans
Let us now indicate the general form of the expression of the derivative
$\left(\tfrac{d^{2n+1}Q}{d\beta^{2n+1}}\right)_0$; first if we find a
term that contains $\tfrac{A_kA_l\ldots A_t}{\Pi k\cdot\Pi l\cdots\Pi t}$
in a factor, it is that we have
\fi
$$
(k+2)+(l+2)+\cdots+(t+2)=2n.
$$
\ifTrans
It remains to determine the coefficient of 
$\tfrac{A_kA_l\ldots}{\Pi k\cdot\Pi l\cdots}$; this coefficient is
composed of different parts joined by the sign of the addition, and we
obtain any one of them in the following way.

Let us write consecutive numbers
\fi
\begin{equation}
\tag{A}\label{eq:9A}
2,\quad 3,\quad 4,\quad 5,\ldots,\quad 2n-1,
\end{equation}
\ifTrans
then suppose that we put in a parenthesis $k$ consecutive of these
numbers, then in a second parenthesis $l$ other consecutive numbers taken
from the numbers \eqref{eq:9A}; then put in a third parenthesis $m$ other
consecutive numbers, and so on. Let us imagine that these parentheses are
separated by at least two of the numbers \eqref{eq:9A}, and always by an
even number of numbers \eqref{eq:9A}; finally let us make this
restriction again, that if the first parenthesis on the left does not
start with the factor $2$, there is before it an even number of numbers 
\eqref{eq:9A}. We will have any part of the coefficient sought, by
multiplying between them the products of the numbers contained in each
parenthesis.

Let us go to the development of $Q_2$. The first derivative of $Q_2$, is
zero for $\beta=0$, and we recognize that it is the same for all the
derivatives of odd order, and by representing by $D$ the value of $Q_2$
for $\beta=0$, we have for even order derivatives
\fi
$$
\def\arraystretch{1.6}
\begin{array}{p{2.5em}p{0.2em}p{22.2em}}
$\left(\tfrac{d^2Q}{d\beta^2}\right)_0$&$=$&$A_0D,\quad
\left(\tfrac{d^4Q}{d\beta^4}\right)_0=
D\left(A_0\tfrac{1\cdot2}{1\cdot2}A_2\right),$\\
$\left(\tfrac{d^6Q}{d\beta^6}\right)_0$&$=$&$D\left(A_0^2+
\tfrac{1\cdot2+3\cdot4}{1\cdot2}A_0A_2+A_4\right),$\\
$\left(\tfrac{d^8Q}{d\beta^8}\right)_0$&$=$&$D\left(A_0^3+
\tfrac{1\cdot2+3\cdot4+5\cdot6}{1\cdot2}A_0^2A_2+
\tfrac{1\cdot2\cdot3\cdot4+3\cdot4\cdot5\cdot6}{1\cdot2\cdot3\cdot4}A_0A_4+\right.$\\
&&$
\qquad\qquad\qquad\qquad\qquad\quad\left.\tfrac{1\cdot2\times 5\cdot6}{1\cdot2\times 1\cdot2}A_2^2+
\tfrac{1\cdot2\cdot3\cdot4\cdot5\cdot6}{1\cdot2\cdot3\cdot4\cdot5\cdot6}A_6\right).$
\end{array}
$$
$$
.\ .\ .\ .\ .\ .\ .\ .\ .\ .\ .\ .\ .\ .\ .\ .\ .\ .\ .\ .\ .\ .\ .\ .\ .\ .\ .\ .\ .\ .\ .\ .\ .\ .\ .\ .\ .\ .\ .\ .\ .\ .\ .\ .
$$
\setcounter{footnote}{0}
\ifTrans
Let us indicate the general form of the derivative 
$\left(\tfrac{d^{2n}Q}{d\beta^{2n}}\right)$. For the term\footnote{Trans.
note: I am not certain this formula is correct, because of the odd
spacing in the original; I have added an ellipsis to attempt to correct it.}
\fi
$$
M\tfrac{A_aA_bA_c\ldots}{\Pi a\cdot\Pi b\cdot \Pi c\cdots},
$$
\ifTrans
$a$, $b$, $c$,\ldots being equal or unequal, we must have
\fi
$$
(a+2)+(b+2)+\ldots=2n;
$$
\ifTrans
it only remains to give the value of the coefficient $M$. To this end,
let's write the numbers
\fi
\begin{equation}
\tag{B}\label{eq:9B}
1,\quad 2,\quad 3,\quad 4,\ldots,\quad 2n-2;
\end{equation}
\ifTrans
this coefficient will be made up of several different parts, each of
which will be obtained as follows. Let us put in a parenthesis $a$
consecutive numbers \eqref{eq:9B}, then in a second parenthesis with $b$
other consecutive numbers, and so on. Let us further imagine that these
parentheses are separated at least by two numbers, and always by an even
number of the numbers \eqref{eq:9B}; finally add this restriction, that
if the first parenthesis does not start with the number $1$, there is
before it an even number of numbers \eqref{eq:9B}. We will have the
sought-after part of the coefficient $M$, by multiplying between them the
products of the numbers contained in each parenthesis.

{\bf 25}. We can develop $P_1$ and $P_2$ in the same way as the functions
$Q_1$ and $Q_2$.

The series obtained for $Q_1$ and $Q_2$ can be written by letting
\fi
$$
M=R-2h^2,
$$
$$
\def\arraystretch{1.6}
\begin{array}{p{0.5em}p{0.2em}p{18em}}
$Q_1$&$=$&$B\left[\beta+M\tfrac{\beta^3}{1\cdot2\cdot3}+
(M^2-24h^2)\tfrac{\beta^5}{1\cdot2\cdot3\cdot4\cdot5}\right.$\\
&&$\qquad\left.+(M^3-104h^2M-160h^2)\tfrac{\beta^7}{1\cdot2\cdot3\cdot4\cdot5\cdot6\cdot7}+\ldots\right],$\\
$Q_2$&$=$&$D\left[1+M\tfrac{\beta^2}{1\cdot2}+
(M^2-8h^2)\tfrac{\beta^4}{1\cdot2\cdot3\cdot4}\right.$\\
&&$\qquad\left.+(M^3-56hM-32h^3)\tfrac{\beta^6}{1\cdot2\cdot3\cdot4\cdot5\cdot6}+\ldots\right].$\\
\end{array}
$$
\ifTrans
By changing $\beta$ to $\alpha i$, we have
\fi
$$
\def\arraystretch{1.6}
\begin{array}{p{0.5em}p{0.2em}p{18em}}
$P_1$&$=$&$B'\left[\alpha-M\tfrac{\alpha^3}{1\cdot2\cdot3}+
(M^2-24h^2)\tfrac{\alpha^5}{1\cdot2\cdot3\cdot4\cdot5}\right.$\\
&&$\qquad\left.-(M^3-104h^2M-160h^2)\tfrac{\alpha^7}{1\cdot2\cdots7}+\ldots\right],$\\
$P_2$&$=$&$D'\left[1-M\tfrac{\alpha^2}{1\cdot2}+
(M^2-8h^2)\tfrac{\alpha^4}{1\cdot2\cdot3\cdot4}\right.$\\
&&$\qquad\left.-(M^3+56hM-32h^3)\tfrac{\alpha^6}{1\cdot2\cdots6}+\ldots\right].$\\
\end{array}
$$

\ifTrans
The value of $R$ or $M$ must be chosen so that $P_1$ and $P_2$ have 
$2\pi$ for period; therefore the values of these two expressions must
remain constant, when we replace $\alpha$ with $\alpha+2\pi$; a very
simple way to determine $M$ is to notice that $P_1$ must become zero for 
$\alpha=\pi$ as for $\alpha=0$, and that $P_2$ must remain the same for
these two values of $\alpha$; we thus have one of the two equations
\fi
\begin{equation}
\tag{$a$}\label{eq:9a}
\left\{
\begin{array}{l}
\pi-M\tfrac{\pi^3}{1\cdot2\cdot3}+(M^2-24h^2)\tfrac{\pi^5}{1\cdot2\cdot3\cdot4\cdot5}\\
\ \ \,-(M^3-104h^2M-160h^2)\tfrac{\pi^7}{1\cdot2\cdot3\cdot4\cdot5\cdot6\cdot7}+\ldots=0,
\end{array}
\right.
\end{equation}
\begin{equation}
\tag{$b$}\label{eq:9b}
\left\{
\begin{array}{l}
1-M\tfrac{\pi^2}{1\cdot2}+(M^2-8h^2)\tfrac{\pi^4}{1\cdot2\cdot3\cdot4}\\
\ \ \hspace{0.5pt}-\,(M^3-56hM-32h^3)\tfrac{\pi^6}{1\cdot2\cdot3\cdot4\cdot5\cdot6}+\ldots=1.\qquad
\end{array}
\right.
\end{equation}

\ifTrans
Suppose for example that it is a vibratory movement of the first kind
given by the formula
\fi
$$
w=P_1Q_1\sin2\lambda mt.
$$

\ifTrans
Let $\beta=\vartheta$ be the equation of the contour which is fixed, $M$
and $h$ will be provided by \eqref{eq:9a} and the equation
\fi
\begin{equation}
\tag{$c$}\label{eq:9c}
\left\{
\begin{array}{l}
\vartheta-M\tfrac{\vartheta^3}{1\cdot2\cdot3}+(M^2-24h^2)\tfrac{\vartheta^5}{1\cdot2\cdot3\cdot4\cdot5}\\
\ \ \,+\,(M^3-104h^2M-160h^2)\tfrac{\vartheta^7}{1\cdot2\cdot3\cdot4\cdot5\cdot6\cdot7}+\ldots=0,
\end{array}
\right.
\end{equation}
\ifTrans
If we especially have in mind the comparison of theory with experience,
we can proceed as follows. After having produced experimentally a
vibratory state of the membrane, one will note the pitch of the sound,
and, consequently, the value of $\lambda=\tfrac{h}{c}$. Then the equation
\eqref{eq:9a} will contain only the unknown $M$, and it will remain to
verify that $M$ and $h=\lambda c$ satisfy \eqref{eq:9c}.

The previous expressions of $P_2$ and $Q_2$ allow to recognize that the
parts of the major axis located between the foci and the neighboring
vertices produce vibrations of maximum amplitude and the part located
between the foci of vibrations of minimum amplitude. Indeed, the value of
$M$ entering it is positive; because we have demonstrated (n$^{\circ}${\bf 21}) 
that $R$ is $>24h^2$. Let us therefore take on the long axis
between the focus and the neighboring vertex a point $n$ for which 
$\alpha$ is zero; consider a very similar point $n'$ on the confocal
ellipse which passes through $n$; $\beta$ is the same for these two
points and $\alpha$ is zero for $n$, very small for $n'$; therefore the
vibratory displacement is greater for $n$ than for $n'$.

Let us take a point $m$ on the line $FF'$ which joins the foci, and also
a point $m'$ very close on the confocal hyperbola which passes by $m$; 
$\alpha$ is the same for $m$ and $m'$, $\beta$ is zero for $m$, very
small for $m'$; therefore the magnitude of the vibration is smaller in 
$m$ than in $m'$.

We have new expressions of $P_1$ and $P_2$ by changing in the previous
values of $P_1$ or $P_2$ according to the parity of $g$ (n$^{\circ}${\bf 14}) 
$\alpha$ to $\tfrac{\pi}{2}-\alpha$, $h^2$ to $-h^2$
(consequently $M$ to $R+2h^2$), and we easily conclude that the minor
axis of the membrane is immobile or at maximum vibration.

\begin{center}
\emph{Annular membrane.}
\end{center}

{\bf 26}. In the two equations
\fi
$$
\rho=c\tfrac{e^{\beta}\,+\,e^{-\beta}}{2},
$$
$$
\tfrac{d^2Q}{d\beta^2}-[R-2h^2E(2\beta)]Q=0,
$$
\ifTrans
make
\fi
$$
\beta=\varepsilon-l\tfrac{c}{2a},
$$
\ifTrans
and we will have the two other
equations
\fi
$$
\rho=a(e^{\varepsilon}+qe^{-\varepsilon}),
$$
$$
\tfrac{d^2Q}{d\varepsilon^2}-[R-f^2(e^{2\varepsilon}+qe^{-2\varepsilon})]Q=0,
$$
\ifTrans
by letting
\fi
$$
2\lambda a=f,\quad \tfrac{c^2}{4a^2}=q,
$$
\ifTrans
and the last two have this advantage over the first that they apply
immediately to the circle by making $q=0$.

Suppose that $Q$ is zero on the ellipse $\rho=\varrho$, let us determine
$a$ so that $\varepsilon$ is zero on this ellipse, we will have to let
\fi
$$
a+\tfrac{c^2}{4a}=\varrho,
$$
\ifTrans
from which
\fi
$$
a=\tfrac{\varrho}{2}\pm\tfrac{\sqrt{\varrho\,-\,c^2}}{2}.
$$ 

\ifTrans
Imagine a ring-shaped membrane whose two fixed edges are confocal
ellipses; if we denote by $\rho=\varrho$ the inner contour, the function
$Q$ develops as follows:
\fi
$$
Q=\varepsilon\left(\tfrac{dQ}{d\varepsilon}\right)_0+
\tfrac{\varepsilon^2}{1\cdot2}\left(\tfrac{d^2Q}{d\varepsilon^2}\right)_0+
\tfrac{\varepsilon^3}{1\cdot2\cdot3}\left(\tfrac{d^3Q}{d\varepsilon^3}\right)_0+\ldots,
$$
\ifTrans
and it remains to determine the expressions of these derivatives.
Although we have here derivatives of even order and others of odd order,
we can form them identically as the derivatives of $Q_1$ with respect to
$\beta$ for $\beta=0$. However, let us content ourselves with writing the
first coefficients of this series in the most convenient form for
numerical calculation:
\fi
$$
\def\arraystretch{1.6}
\begin{array}{p{3.5em}p{0.2em}p{26em}}
$\quad\ \left(\tfrac{dQ}{d\varepsilon}\right)_0$&$=$&$B,\quad
\left(\tfrac{d^2Q}{d\varepsilon^2}\right)_0=0,\quad
\tfrac{1}{B}\left(\tfrac{d^3Q}{d\varepsilon^3}\right)_0=-f^2(1+q)+R,$\\
$\tfrac{1}{B}\left(\tfrac{d^4Q}{d\varepsilon^4}\right)_0$&$=$&$-4f^2(1-q),$\\
$\tfrac{1}{B}\left(\tfrac{d^5Q}{d\varepsilon^5}\right)_0$&$=$&$f^4(1+q)^2-
2f^2(R+6)(1+q)+R^2,$\\
$\tfrac{1}{B}\left(\tfrac{d^6Q}{d\varepsilon^6}\right)_0$&$=$&$
12f^4(1+q)^2-4f^2(1-q)(3R+8),$\\
$\tfrac{1}{B}\left(\tfrac{d^7Q}{d\varepsilon^7}\right)_0$&$=$&$
-f^6(1+q)^3+f^4[3R(1+q)^2+4(23+6q+23q^2)],$\\
&&$-f^2(1+q)(3R^2+52R+80)+R^3,$\\
$\tfrac{1}{B}\left(\tfrac{d^8Q}{d\varepsilon^8}\right)_0$&$=$&$
-24f^6(1-q)(1+q)^2+48f^4(1-q^2)(R+12)$\\
&&$-24f^2(1-q)(R^2+8R+8),$\\
$\tfrac{1}{B}\left(\tfrac{d^9Q}{d\varepsilon^9}\right)_0$&$=$&$
f^8(1+q)^4-4f^6(1+q)[R(1+q)^2+86-36q+86q^2]$\\
&&$+f^4[6R^2(1+q)^2+32R(15+4q+15q^2)$\\
&&$\qquad\qquad\qquad\qquad\qquad\qquad\,+\,16(201+10q+201q^2)]$\\
&&$-4f^2(1+q)(R^3+34R^2+160R+112)+R^4,$
\end{array}
\vspace{-7pt}
$$
$$
\def\arraystretch{1.6}
\begin{array}{p{3.5em}p{0.2em}p{26em}}
$\tfrac{1}{B}\left(\tfrac{d^{10}Q}{d\varepsilon^{10}}\right)_0$&$=$&$
40f^8(1-q)(1+q)^3$\\
&&$-f^6(1-q)(1+q)^2\left[120R+3200+640\left(\tfrac{1-q}{1+q}\right)^2\right]$\\
&&$+f^4(1-q^2)(120R^2+3840R+16704)$\\
&&$-f^2(1-q)(40R^3+640R^2+1984R+1024),$\\
$\tfrac{1}{B}\left(\tfrac{d^{11}Q}{d\varepsilon^{11}}\right)_0$&$=$&$
-f^{10}(1+q)^5+5f^8(1+q)^2[R(1+q)^2+184-144q+184q^2]$\\
&&$-f^6(1+q)[10R^2(1+q)^2+40R(53-22q+53q^2)$\\
&&$\qquad\qquad\qquad\qquad\qquad\,+\,36912-29216q+36912q^2]$\\
&&$+f^4[10R(1+q)^2+R^2(1480+400q+1480q^2)$\\
&&$\qquad+R(26896+1440q+26896q^2)$\\
&&$\qquad\qquad\quad\qquad\quad\,+\,82624+896q+82624q^2]$\\
&&$-f^2(1+q)(5R^4+280R^3+2656R^2+5824R+2304)+R^5.$
\end{array}
$$

\ifTrans
According to whether the constant $R$ is relative to a function $P_1$ or
$P_2$, we have for the vibratory movement of the annular membrane
\fi
$$
w=P_1Q\sin2\lambda mt
$$
\ifTrans
or
\fi
$$
w=P_2Q\sin2\lambda mt,
$$
\ifTrans
by giving $Q$ the value we just calculated; then we finish determining
the movement by calculating the quantity $\lambda$ according to the
condition that $Q$ is zero on the external contour $\rho=A$.

{\bf 27}. Thus there are two kinds of solutions: in one, the portions of
the major axis located between the two fixed contours are nodes, and in
the other, are vibration bellies, and the nodal lines are still ellipses
and portions of confocal hyperbolas.

When it is a full membrane, we have two distinct solutions
\fi
$$
w=P_1Q_1\sin2\lambda mt,\quad w=P_2Q_2\sin2\lambda mt,
$$
\ifTrans
and $Q_1$ and $Q_2$ have two distinct forms, like $P_1$ and $P_2$. On the
contrary, when the membrane is annular, the two functions $Q$ which
associate with $P_1$ and $P_2$ only differ by the constant $R$ which
enters it. It follows from this that if the integer $g$, which designates
the number of hyperbolic nodal lines, is large enough, and especially if
at the same time the eccentricity is small enough, the constant $R$ will
differ very little for an identical value of $g$ in the two functions
$P_1$ and $P_2$, and the two functions $Q$ associated with them will be
almost identical. The two corresponding vibratory states will therefore
render almost the same sound, and, according to an experiment, they will
be superimposed, and the resulting state will be represented by the
formula
\fi
$$
W=(AP_1+BP_2)Q\sin2\lambda mt;
$$
\ifTrans
then the hyperbolic nodal lines will be given by the equation
\fi
$$
AP_1+BP_2=0
$$
\ifTrans
and will still be $g$ in number. But previously the same hyperbola
provided portions of these four branches; now a hyperbola provides only
two portions of branches which have the same asymptote; because the
function $AP_1+BP_2$ is not symmetrical with respect to the axes of the
ellipse, but if we change $\alpha$ to $\pi+\alpha$, it remains the same
or changes sign depending on whether $g$ is even or odd.
\fi

\ifTrans
We see that then the sound and the elliptical nodal lines remaining
invariable, the position of the hyperbolic lines which depends on the
arbitrary ratio $\tfrac{B}{A}$ can vary although their number $g$ does
not change.

It is obvious that the previous formulas apply to the circular ring by
making $q=0$, $R=g^2$, and they are also suitable for the full elliptical
membrane for movement of the first kind
\fi
$$
w=P_1Q_1\sin2\lambda mt,
$$
\ifTrans
since it suffices to assume that the fixed interior contour is reduced to
the line segment 
of the focal points, which is the limit of the smallest
confocal ellipses. So you have to make $\varrho=c$, $a=\tfrac{c}{2}$, 
$q=1$, and $\varepsilon$ is reduced to $\beta$.

\begin{center}
\emph{Elliptical nodal lines.}
\end{center}

{\bf 28}. Consider the values of $Q_1$ and $Q_2$ given by the equation
\fi
\begin{equation}
\tag{1}\label{eq:101}
\tfrac{d^2Q}{d\beta^2}+[h^2(e^{2\beta}+e^{-2\beta})-R]Q=0,
\end{equation}
\ifTrans
one of which is zero and the other minimum for $\beta=0$.

Suppose that we increase $\lambda$ and consequently $h=\lambda c$, since
$c$ is fixed: $Q_1$ and $Q_2$ vary, but keeping their character at the
limit $\beta=0$. $R$ is a function of $h$, and let us start by assuming
that
\fi
$$
\tfrac{dR}{dh}-4h\quad\rm{est}\quad<0;
$$
\ifTrans
then, as $\lambda$ increases, the coefficient of $Q$ in \eqref{eq:101} will
take larger and larger values; because from the previous inequality we
conclude that the derivative of this coefficient with respect to $h$,
\fi
$$
2h(e^{2\beta}+e^{-2\beta})-\tfrac{dR}{dh},
$$
\ifTrans
is $>0$. So the roots of the equation in $\beta$
\fi
$$
Q(\beta,\lambda)=0
$$
\ifTrans
decrease in size as $\lambda$ increases.
\fi

\ifTrans
Let us denote by $\beta=B$ the parameter of the contour ellipse, on which
$Q$ is zero, the equation
\fi
$$
Q(B,\lambda)=0
$$
\ifTrans
determines the number $\lambda$. Let $\lambda_1$, $\lambda_2$, 
$\lambda_3$,\ldots be these roots in order of increasing size: 
$\lambda_i$ being the $i^{th}$ root, $Q(\beta,\lambda_i)$ is one of
the values of $Q$ from our research, and the equation
\fi
$$
Q(\beta,\lambda_i)=0
$$
\ifTrans
will give, by its roots in $\beta$, the parameters of the elliptical
nodal lines. Now we will prove that this equation has $i-1$ roots, 
$\beta_1$, $\beta_2$,\ldots, $\beta_{i-1}$, less than $B$, and that
consequently the nodal ellipses are equal to $i-1$ in number.

Consider the equation in $\beta$,
\fi
\begin{equation}
\tag{$a$}\label{eq:10a}
Q(\beta,\lambda)=0
\end{equation}
\ifTrans
and represent the curve
\fi
$$
y=Q(\beta,\lambda),
$$
\ifTrans
$\beta$ being taken for abscissa and $y$ for the ordinate, which is zero
or maximum for $\beta=0$. Let $i$ be the number of roots between $0$ and
$B$, which are determined by the points $\beta_1$, $\beta_2$,\ldots, 
$\beta_i$. Let us increase $\lambda$,
\fi
\begin{figure}[h!]
\begin{center}
\includegraphics[width=0.6\textwidth]{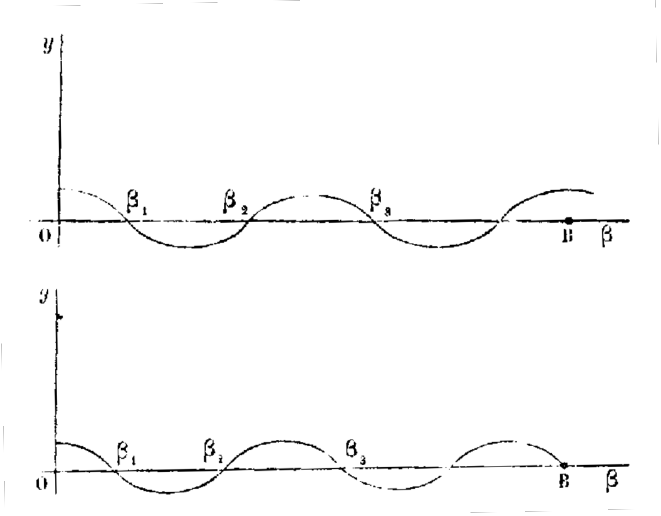}
\end{center}
\begin{caption}\\
\footnotesize Image source: \url{http://sites.mathdoc.fr/JMPA/PDF/JMPA_1868_2_13_A8_0.pdf}
\end{caption}
\end{figure}
\ifTrans
the points $\beta_1$, $\beta_2$,\ldots will approach the origin, the
sinuosities of the curve will decrease in amplitude, and for a value 
$\lambda_{i+1}$ the curve will pass through the point $B$. So we have a
value of $\lambda=\lambda_{i+1}$ which
\fi
\ifTrans
satisfies the equation in $\lambda$
\fi
\begin{equation}
\tag{$b$}\label{eq:10b}
Q(B,\lambda)=0;
\end{equation}
\ifTrans
and if we give a new small increase to $\lambda$, the equation 
\eqref{eq:10a} will have a new root to the left of $B$, and will
therefore have $i+1$ roots between $0$ and $B$.

Let us continue to increase $\lambda$, the points $\beta_1$, $\beta_2$,\ldots
are getting closer to zero again, and, for a value 
$\lambda=\lambda_{i+2}$, the curve will pass again by the point $B$; we
will therefore have another new value of $\lambda$, $\lambda_{i+2}$ which
satisfies \eqref{eq:10b}, and the equation
\fi
$$
Q(\beta,\lambda_{i+2}+\varepsilon)=0
$$
\ifTrans
where $\varepsilon$ is positive and very small, $i+2$ roots between $0$ and
$B$; the equation
\fi
$$
Q(\beta,\lambda_{i+2})=0
$$
\ifTrans
itself will have $i+2$, counting $B$. Now $\lambda_{i+1}$ and 
$\lambda_{i+2}$ are obviously two consecutive roots of the equation 
\eqref{eq:10b} in $\lambda$; we conclude that if $\lambda_{i+1}$ and 
$\lambda_{i+2}$ are two consecutive roots of \eqref{eq:10b}, the equation
in $\beta$
\fi
$$
Q(\beta,\lambda_{i+2})=0
$$
\ifTrans
has one more root than
\fi
$$
Q(\beta,\lambda_{i+1})=0
$$
\ifTrans
between the limits $0$ and $B$. That said, we can easily recognize that 
$Q(\beta,\lambda_1)=0$ has no roots between $0$ and $B$; so 
$Q(\beta,\lambda_2)=0$ has one, $Q(\beta,\lambda_3)=0$ has two, etc., and
in general $Q(\beta,\lambda_i)=0$ has $i-1$.

The number of hyperbolic nodal lines remaining the same, we see that, as
the sound rises, the number of elliptical nodal lines increases. All of
the above is based on the existence of the inequality
\fi
\begin{equation}
\tag{$c$}\label{eq:10c}
\tfrac{dR}{dh}-4h<0.
\end{equation}

\ifTrans
At n$^{\circ}${\bf 8}, we found the equation
\fi
$$
\tfrac{dP}{d\alpha}\delta P-P\delta\tfrac{dP}{d\alpha} = 
\delta h\int_0^{\alpha}P^2\left(\tfrac{dR}{dh}-4h\cos2\alpha\right)d\alpha;
$$
\ifTrans
for $\alpha=0$, the two members are zero; but if we suppose $\alpha$
excessively small, $\tfrac{dR}{dh}-4h\cos2\alpha$ will not change sign
between $0$ and this value of $\alpha$; so the first member will have the
same sign as $\tfrac{dR}{dh}-4h$.

If this quantity can be positive, the expression
\fi
$$
\tfrac{dR}{dh}-2h(e^{2\beta}+e^{-2\beta})
$$
\ifTrans
will also be 
for excessively small values of $\beta$;
therefore by taking $B$ sufficiently small, and, consequently, the 
very-eccentric membrane, the number of elliptical nodal lines would
decrease when, the number of hyperbolic nodal lines remaining the same,
the pitch of the sounds would increase. As this result does not seem
admissible, it seems that the inequality \eqref{eq:10c} must always be
tied-in; however, the above cannot be viewed as a rigorous demonstration.

\newpage
\begin{center}
\emph{Most general vibratory movement of the elliptical membrane.}
\end{center}

{\bf 29}. We have so far dealt only with simple vibratory movements,
which are those which would be most easily produced in experience. We
will now assume that we give at all points of a membrane any initial
velocities, and determine the vibrational state that will result.

But first let us do some thinking about the signs of the coordinates we
use. As we said at n$^{\circ}${\bf 4}, where we let
\fi
\begin{equation}
\tag{$a$}\label{eq:11a}
\left\{
\begin{array}{p{0.1em}p{0.2em}p{10em}}
$x$&$=$&$cE(\beta)\cos\alpha,$\\
$y$&$=$&$c\,\mathcal{E}(\beta)\sin\alpha,$
\end{array}
\right.
\end{equation}
\ifTrans
when using the coordinates $\alpha$ and $\beta$, we can assume that 
$\beta$
\fi
\ifTrans
is essentially positive, and that $\alpha$ is only susceptible to varying
from zero to $2\pi$, or from $-\pi$ to $+\pi$, and despite these
restrictions we can represent by these coordinates any point of the plane.

But it follows from the formulas \eqref{eq:11a}, that if we give a
negative value to $\beta$ instead of giving it to $\alpha$, the point 
$(x,y)$ remains the same, and therefore also that the coordinates 
$(-\alpha,-\beta)$ represent the same point as the coordinates 
$(\alpha,\beta)$. So a formula which gives the vibratory movement of the
membrane must remain invariable when we replace $\alpha$ and $\beta$ with
$-\alpha$ and $-\beta$: this is what we verify indeed on the two simple
solutions that we found
\fi
$$
w=P_1Q_1\sin2\lambda mt,\quad w=P_2Q_2\sin2\lambda mt,
$$
\ifTrans
since $P_1$ and $Q_1$ are odd functions of $\alpha$ and $\beta$, and
since $P_2$ and $Q_2$ are even functions, and we could have associated
the functions $Q_1$ and $Q_2$ to functions $P_1$ and $P_2$, according to
this condition (n$^{\circ}${\bf 15}).

However, it should be noted that these considerations would not apply to
the annular membrane. Indeed, the line drawn between the foci, and which
has the equation $\beta=0$ is no longer located on the surface of the
membrane, and if we have expressed that $Q$ is zero for $\beta=\beta_1$ on the inside contour
and taken $\beta_1$ positive, we can only
give $\beta$ the positive values enclosed between those which suit the
two contours.

Returning to the vibratory movement of the full membrane, suppose that
the initial speed given at each point of the membrane is expressed by the
formula
\fi
$$
\left(\tfrac{dw}{dt}\right)_0=\Phi(\alpha,\beta),
$$
\ifTrans
in which $\Phi(\alpha,\beta)$ is a function which becomes zero on the contour
of the membrane $\beta=\vartheta$, and which, from what we have seen,
remains invariable when we replace $\alpha$ and $\beta$ with $-\alpha$
and $-\beta$. We easily conclude that $\Phi(\alpha,\beta)$ is the sum of
two functions $F_1(\alpha,\beta)$, $F_2(\alpha,\beta)$, which, ordered
with respect to increasing powers of $\alpha$ and $\beta$, are: one of
the form  
\fi
$$
F_2=a+A\alpha^2+B\beta^2+C\alpha^4+D\alpha^2\beta^2+
E\beta^4+F\alpha^6+G\alpha^4\beta^2+\ldots.
$$
\ifTrans
even in $\alpha$ and $\beta$; and the other of the form
\fi
$$
F_1=A'\alpha\beta+B'\alpha^3\beta+C'\alpha\beta^3+D'\alpha^5\beta+
E'\alpha^3\beta^3+F'\alpha\beta^5+G'\alpha^7\beta+\ldots,
$$
\ifTrans
odd in $\alpha$ and odd in $\beta$, but even with respect to their set.

After having let
\fi
$$
\Phi(\alpha,\beta)=F_1(\alpha,\beta)+F_2(\alpha,\beta),
$$
\ifTrans
let us look at the resulting vibratory movement as the sum of an infinity
of simple vibratory movements, the amplitude of which we will propose to
determine. Each simple movement of the first or second kind given by the
formulas
\fi
$$
w=aP_1Q_1\sin2\lambda mt,\quad w=bP_2Q_2\sin2\lambda mt
$$
\ifTrans
depends first on an integer $g$, and, this number $g$ once designated,
this movement can vary in an infinite number of ways by the number 
$\lambda$, which is susceptible to increasing values $\lambda_1$, 
$\lambda_2$,\ldots, $\lambda_i$,\ldots, and we will assign them a second
index which recalls the number $g$, and we will replace the two previous
formulas by the following two:
\fi
$$
\def\arraystretch{1.6}
\begin{array}{p{0.2em}p{0.2em}p{14em}}
$w$&$=$&$aP_1(g,\lambda_i^g)Q_1(g,\lambda_i^g)\sin2\lambda_i^gmt,$\\
$w$&$=$&$bP_2(g,\lambda'_i{}^g)Q_1(g,\lambda'_i{}^g)\sin2\lambda'_i{}^gmt,$
\end{array}
$$
\ifTrans
Then considering a vibratory state composed of an infinity of simple states, we will have
\fi
$$
\begin{array}{p{0.4em}p{20em}}
$w$&$=\sum a_{g,\lambda_i}P_1(g,\lambda_i^g)Q_1(g,\lambda_i^g)\sin2\lambda_i^gmt$\\
&$+\,\sum b_{g,\lambda'_i}P_2(g,\lambda'_i{}^g)Q_2(g,\lambda'_i{}^g)\sin2\lambda'_i{}^gmt$
\end{array}
$$
\ifTrans
and we draw for the initial speed
\fi
$$
\begin{array}{p{2em}p{20em}}
$\left(\tfrac{dw}{dt}\right)_0$&$=2m\sum\lambda_i^g a_{g,\lambda_i}P_1(g,\lambda_i^g)Q_1(g,\lambda_i^g)$\\
&$+\ 2m\sum\lambda'_i{}^g b_{g,\lambda'_i}P_2(g,\lambda'_i{}^g)Q_2(g,\lambda'_i{}^g),$
\end{array}
$$
\ifTrans
expression which must be identified with $\Phi(\alpha,\beta)$; but we
will decompose this equality into the following two
\fi
\begin{equation}
\tag{1}\label{eq:111}
F_1(\alpha,\beta)=2m\sum\lambda_i^g a_{g,\lambda_i}P(g,\lambda_i^g)Q(g,\lambda_i^g),
\end{equation}
\begin{equation}
\tag{2}\label{eq:112}
\ \ F_2(\alpha,\beta)=2m\sum\lambda'_i{}^g b_{g,\lambda'_i}P(g,\lambda'_i{}^g)Q(g,\lambda'_i{}^g).
\end{equation}

\ifTrans
Now consider the four equations
\fi
\begin{equation}
\tag{$b$}\label{eq:11b}
\left\{
\def\arraystretch{1.6}
\begin{array}{c}
\tfrac{d^2Q}{d\beta^2}\,-\,[R(g,\lambda c)\,-\,2\lambda^2c^2E(2\beta)]\,Q\ =0,\\
\tfrac{d^2Q'}{d\beta^2}-[R(g',\lambda' c)-2\lambda'^2c^2E(2\beta)]Q'=0;
\end{array}
\right.
\end{equation}
\begin{equation}
\tag{$c$}\label{eq:11c}
\left\{
\def\arraystretch{1.6}
\begin{array}{c}
\tfrac{d^2P}{d\alpha^2}\,+\,[R(g,\lambda c)\,-\,2\lambda^2c^2\cos2\alpha]\,P\ \,=0,\\
\tfrac{d^2P'}{d\alpha^2}+\,[R(g',\lambda' c)-2\lambda'^2c^2\cos2\alpha]P'=0.
\end{array}
\right.
\end{equation}

\ifTrans
By subtracting the two equations \eqref{eq:11b} multiplied by $Q'$ and 
$Q$, we have
\fi
$$
0=Q'\tfrac{d^2Q}{d\beta^2}-Q\tfrac{d^2Q'}{d\beta^2}+
[2(\lambda^2-\lambda')c^2E(2\beta)-(R-R')]QQ';
$$
\ifTrans
we integrate from $\beta=0$ to $\beta=\vartheta$, the parameter of the
contour, and we will have
\fi
\begin{multline}
\nonumber
0=\left(Q'\tfrac{dQ}{d\beta}-Q\tfrac{dQ'}{d\beta}\right)_{\vartheta}-
\left(Q'\tfrac{dQ}{d\beta}-Q\tfrac{dQ'}{d\beta}\right)_{0}\\
+2(\lambda^2-\lambda'^2)c^2\displaystyle\int_0^{\vartheta}E(2\beta)QQ'\,d\beta-(R-R')\displaystyle\int_0^{\vartheta}QQ'\,d\beta.\qquad
\end{multline}
\ifTrans
The first term is zero because $Q$ and $Q'$ are zero for 
$\beta=\vartheta$; then, if $Q$ and $Q'$ have the character of $Q_1$,
they are zero for $\beta=0$, and if they both have the character of 
$Q_2$, their derivatives are zero for $\beta=0$; so the second term is
also zero. We also find
\fi
\begin{multline}
\nonumber
0=\left(P'\tfrac{dP}{d\alpha}-P\tfrac{dP'}{d\alpha}\right)_{0}^{2\pi}+
2(\lambda^2-\lambda'^2)c^2\displaystyle\int_0^{2\pi}PP'\cos2\alpha\,d\alpha\\
\quad-(R-R')\displaystyle\int_0^{2\pi}PP'\,d\alpha.\qquad\qquad\qquad\qquad\qquad\qquad
\end{multline}
\ifTrans
whose first part is zero, because $P$ and $P'$ are periodic functions. So
we have the two equalities
\fi
\begin{equation}
\tag{$d$}\label{eq:11d}
\left\{
\begin{array}{c}
(R-R')\displaystyle\int_0^{\vartheta}QQ'\,d\beta=2(\lambda^2-\lambda'^2)c^2\displaystyle\int_0^{\vartheta}QQ'E(2\beta)\,d\beta,\\
2(\lambda^2-\lambda'^2)c^2\displaystyle\int_0^{2\pi}PP'\cos2\alpha\,d\alpha=(R-R')\displaystyle\int_0^{2\pi}PP'\,d\alpha.
\end{array}
\right.
\end{equation}

\ifTrans
Multiply these equalities member-to-member, and, dividing by
\fi
$$
2(R-R')(\lambda^2-\lambda'^2)c^2
$$
\ifTrans
we obtain
\fi
\begin{equation}
\tag{$e$}\label{eq:11e}
\int_0^{\vartheta}\int_0^{2\pi}[E(2\beta)-\cos2\alpha]PP'QQ'\,d\beta d\alpha=0.
\end{equation}
\ifTrans
This equality is no longer demonstrated if $\lambda=\lambda'$ or if 
$R=R'$; it is however still exact, because, if $\lambda=\lambda'$, we
will deduce from the equations \eqref{eq:11d}
\fi
$$
\int_0^{\vartheta}QQ'\,d\beta=0,\quad
\int_0^{2\pi}PP'\,d\alpha=0;
$$
\ifTrans
therefore the two parts of the integral \eqref{eq:11e} are zero. If 
$R'=R$, we still see that the two equalities \eqref{eq:11d} entail 
\eqref{eq:11e}. 

We multiply the two members of equality \eqref{eq:111} by
\fi
$$
P_1(g,\lambda_i^g)Q_1(g,\lambda_i^g)[E(2\beta)-\cos2\alpha]\,d\alpha d\beta,
$$
\ifTrans
and integrate, with respect to $\alpha$, from $0$ to $2\pi$, and, with
respect to $\beta$, from $0$ to $\vartheta$: all the terms will disappear
in the second member according to \eqref{eq:11e}, except that which has
the coefficient $a_{g,\lambda_i}$ which is determined. We also have 
$b_{g,\lambda_i}$ by means of \eqref{eq:112}.

We have a similar calculation for the annular membrane fixed between two
confocal ellipses; it would be superfluous to insist on it, and even we
did the previous calculation only because it required considerations
relating to the signs of $\alpha$ and $\beta$ that are useful to notice.

\fi
\ifTrans
To return to these signs, let us still imagine that one has to seek the
movement of an elliptical membrane from which one removes the two
portions cut by a confocal hyperbola, and suppose all the contour is
fixed. We will again have a simple vibratory movement represented by the
formula
\fi
$$
w=PQ\sin2\lambda mt;
$$
\ifTrans
but $P$ is no longer a periodic function of $\alpha$. We must vary 
$\alpha$ only between the limits $\alpha_1$, and $\pi-\alpha_1$, relative
to the hyperbola of the contour, and we will vary $\beta$ between the two
limits $-\vartheta$ and $+\vartheta$ relative to the two elliptical arcs
of the periphery of the membrane.
\fi


\end{document}